\documentclass[12pt]{amsart}

\usepackage{amsthm, amsmath, amssymb, comment}
\usepackage{subfig}
\usepackage{graphicx}
\usepackage{subfig}
\usepackage{color}

\newtheorem{theorem}{Theorem}
\newtheorem{lemma}{Lemma}
\newtheorem{proposition}{Proposition}
\newtheorem{corollary}{Corollary}
\theoremstyle{definition}
\newtheorem{definition}{Definition}
\newtheorem{example}{Example}
\theoremstyle{remark}
\newtheorem{remark}{Remark}

\date{}

\author{Avetik Arakelyan}
\address{Institute of Mathematics, National Academy of Sciences of Armenia, 0019 Yerevan, Armenia}
\email{arakelyanavetik@gmail.com}

\author{Rafayel Barkhudaryan}
\address{Institute of Mathematics, National Academy of Sciences of Armenia and American University of Armenia, 0019 Yerevan, Armenia}

\email{rafayel@instmath.sci.am}

\title[A numerical approach for a general class of reaction-diffusion systems ]{A numerical approach for a general class of the spatial segregation of reaction-diffusion systems arising in population dynamics}


\begin{document}

\begin{abstract}
In the current work we consider the numerical solutions of equations of stationary states for a general class of  the spatial segregation of reaction-diffusion systems with $m\geq 2$ population densities. We introduce a discrete multi-phase minimization problem related to the segregation problem, which allows to prove the existence and uniqueness of the corresponding finite difference scheme. Based on that scheme, we suggest an iterative algorithm and show its consistency  and stability. For the special case $m=2,$ we show that the problem gives rise to the generalized version of the so-called two-phase obstacle problem. In this particular case we introduce the notion of viscosity solutions and  prove convergence of the difference scheme to the unique viscosity solution. At the end of the paper we present computational tests,  for different internal dynamics, and discuss numerical results.
\end{abstract}

\maketitle

\section{Introduction and known results}

\subsection{The  setting of the problem }
In recent years there have been intense studies of spatial segregation for reaction-diffusion systems. The existence of spatially inhomogeneous solutions for competition models of Lotka-Volterra type in the case of two and more competing densities  have been considered in
\cite{MR2146353,MR2151234,MR2283921,MR2300320,MR1459414, MR1900331, MR2417905}.
 The aim of this paper is to study the numerical solutions for a certain class of the spatial segregation of reaction-diffusion system with $m$ population densities. 

Let $\Omega  \subset \mathbb{R}^n, n\geq 2$ be a connected and  bounded domain with smooth boundary and  $m$ be a fixed integer.
We consider the steady-states of $m$ competing species coexisting in  the same area $\Omega$.
Let $u_{i}$ denotes  the population density of the  $i^\textrm{th}$ component with the internal dynamic  prescribed by $F_{i}$.

We call the $m$-tuple $U=(u_1,\cdots,u_m)\in (H^{1}(\Omega))^{m},$ \emph{a segregated state} if
\[
u_{i}(x) \cdot  u_{j}(x)=0,\  \text{a.e. } \text{ for  } \quad i\neq j,  \ x\in \Omega.
\]
The problem amounts to
\begin{equation}\label{main_problem}
\text{  Minimize  }  E(u_1, \cdots, u_m)=\int_{\Omega}  \sum_{i=1}^{m} \left( \frac{1}{2}| \nabla u_{i}(x)|^{2}+F_i(x,u_i(x)) \right) dx
\end{equation}
  over the set
  $$S={\{(u_1,\dots,u_{m})\in (H^{1}   (\Omega))^{m} :u_{i}\geq0, u_{i} \cdot u_{j}=0, u_{i}=\phi_{i} \quad \text {on} \quad \partial  \Omega}\},$$
where $\phi_{i} \in H^{\frac{1}{2}}(\partial \Omega),$\; $\phi_{i}  \cdot \phi_{j}=0,$ for $i\neq j$ and $ \phi_{i}\geq 0$ on the boundary $\partial \Omega$.

We assume that
\[F_i(x,s)=\int_0^s f_i(x,v)dv,\]
 where $f_i(x,s):\Omega\times\mathbb{R}^+\to\mathbb{R}$ is Lipschitz continuous in $s,$ uniformly continuous in $x$ and $f_i(x,0)\equiv0$.
%In the sequel, we assume that the functional  \eqref{main_problem} is coercive, which will be needed to provide the existence of minimizers. In order to have  coercivity for the functional, for instance following \cite{ MR2151234}, one can assume that
%there exists $b_i(x)\in L_\infty(\Omega)$ such that
%\[
%|f_i(x,s)|\leq b_i(x)\cdot s, \;\;\forall x\in\Omega,\;\; s\geq s_0 >> 1,
%\]
%and
%\[
%\int_\Omega(|\nabla w(x)|^2 - b_i(x)w^2(x))dx>0, \quad\forall w \in H_0^1(\Omega)\setminus\{0\}.
%\]

\begin{remark}
Functions $f_i$'s are defined only for non negative values of s (recall that our densities $u_i$'s are assumed non negative); thus we can arbitrarily define such functions on the negative semiaxis. For the sake of convenience, when $s\leq 0$,  we will let $f_i(x,s)=-f_i(x,-s)$. This extension  preserves the continuity due to the conditions on $f_i$ defined above. In the same way, each $F_i$ is extended as an even function.
\end{remark}
\begin{remark}
We emphasize that for the case $f_i(x,s)=f_i(x),$ the  assumption is that  for all $i$ the functions $f_i$ are nonnegative and  uniformly continuous in $x$.  Also for simplicity, throughout the paper we shall call both $F_i$ and $f_i$  as internal dynamics.

%And throughout the paper we will  consider the case $f_i(x,s)$ depends on the variable $s,$ keeping in mind that the same results will take place also when  $f_i(x,s)$  doesn't depend on $s$. 
\end{remark}

\begin{remark}
We would like to point out that the only difference between our minimization problem \eqref{main_problem} and the problem discussed by Conti,Terrachini and Verzini \cite{MR2151234}, is  the sign in front of the internal dynamics $F_i$.  In our case, the plus sign of $F_i$ allows to get rid of some additional  conditions, which are imposed in \cite[Section $2$]{MR2151234}. Those conditions are important to provide coercivity of a minimizing functional in \cite{MR2151234}. But in our case the above given conditions together with convexity assumption on $F_i(x,s),$ with respect to the variable $s$ are enough to conclude $F_i(x,u_i(x))\geq 0,$ which in turn implies coercivity of a functional \eqref{main_problem}.
\end{remark}

In order to speak on  the local properties of the population densities, let us introduce the notion of multiplicity of a point in $\Omega$.
\begin{definition}
The multiplicity of the point $x \in \overline{\Omega}$ is defined by:
\begin{equation*}
m(x)=\text{card} \left\{i: measure (\Omega_{i} \cap B(x,r))>0, \forall r>0\right\},
\end{equation*}
where  $\Omega_i=\{u_i>0\}$.
\end{definition}

For the local properties of $u_i$ the same results as in \cite{MR2151234} with  the opposite sign in front of the internal dynamics $f_i$ hold. Below, for the sake of clarity, we  write down those results from \cite{MR2151234} with appropriate changes.
\begin{lemma}(Proposition $6.3$ in \cite{ MR2151234})
Assume that $x_{0} \in \Omega,$ then the following holds:\\
\begin{enumerate}
\item[1)] If $m(x_0)=0,$ then there exists $r>0$ such that for every $i=1,\cdots m$;
\[
 u_i\equiv0\;\; \mbox{on}\;\;  B(x_0,r).
\]
\item[2)] If $m(x_0)=1,$ then there are $i$  and $r>0$ such that in $ B(x_0,r)$
\[
\Delta u_i=f_i(x,u_i), \quad \ \quad u_j\equiv0 \quad \text{for  } j\neq i.
\]
\item[3)] If $m(x_0)=2,$ then there are $i, j $ and  $r>0$ such that for every $k$  and $k\neq i,j$, we have $ u_k\equiv 0$ and in $B(x_0,r)$
\[
\Delta(u_i-u_j) =f_{i}(x,(u_i-u_j))\chi_{\{u_i>u_j\}} -   f_{j}(x,-(u_i-u_j))\chi_{\{u_i<u_j\}}.
\]
\end{enumerate}
\end{lemma}

Next, we state the following uniqueness Theorem due to  Conti, Terrachini and Verzini.
\begin{theorem}[Theorem $4.2$ in \cite{ MR2151234}]
Let the functional in minimization problem \eqref{main_problem} be coercive and moreover each $F_i(x,s)$ be convex in the variable $s,$ for all $x\in\Omega$. Then, the problem  \eqref{main_problem} has a unique minimizer.
\end{theorem}
This theorem will play a crucial role in studying the difference scheme, especially for the case $m=2$ where we will reformulate it as a  generalized two-phase obstacle problem. Note that in this case, the problem will be reduced to:

\begin{equation}\label{m=2}
\text{  Minimize  }  E(u_1, u_2)=\int_{\Omega}  \sum_{i=1}^{2} \left( \frac{1}{2}| \nabla u_{i}(x)|^{2}+F_i(x,u_i(x))\right) dx,
\end{equation}
  over the set
  $$S={\{(u_1,u_{2})\in (H^{1}   (\Omega))^{2} :u_{i}\geq0, u_{1} \cdot u_{2}=0, u_{i}=\phi_{i} \quad \text {on} \quad \partial  \Omega}\}.$$
Here $\phi_{i} \in H^{\frac{1}{2}}(\partial \Omega)$ with property $\phi_{1}  \cdot \phi_{2}=0, \, \phi_{i}\geq 0$ on the boundary $\partial \Omega$. %Also  we assume that $f_i$ is uniformly  continuous   and    $f_{i}(x) \geq 0$.

\subsection{Known results}
In last years there has been much interest given to study the numerical approximation of reaction-diffusion type equations. For instance,  the equations arising  in the study of population ecology when high  competitive interactions between different species occurs.

We refer the reader to \cite{MR2151234,MR2363653, MR2079274, MR1303035, MR1687440, MR1900331, MR1459414} for models involving Dirichlet boundary data. A complete analysis of the stationary case has been studied  in \cite{MR2151234}. Also  numerical simulation   for the spatial segregation limit of two diffusive Lotka-Volterra models in presence of strong competition and inhomogeneous Dirichlet boundary conditions is provided in \cite{MR2459673}.

 In the work \cite{MR2563520} Bozorgnia proposed two numerical algorithms for the problem \eqref{main_problem} with the internal dynamics $F_i\equiv 0$. The finite element approximation is based on the local properties of the solution. In this case the author was able to provide the convergence of the method. The second approach is a finite difference method, but lack of its analysis in \cite{MR2563520}. This finite difference method has been generalized in \cite{Mywork} for the case of non-negative $f_i$. In \cite{Mywork} the authors present a numerical consistent variational system with strong interaction, and provide disjointness condition of populations during the iteration of the scheme. In this case the proposed algorithm is lack of deep analysis,  especially for the case of three and more competing populations.   
 
 The present work concerns to close that gap and provides theoretical results for finite difference scheme, with  $m\geq 2$ competing populations and general  internal dynamics $F_i,$  satisfying certain suitable conditions. We  introduce the discrete analogue of minimization problem and prove the existence and uniqueness of the difference scheme. Moreover for the special case $m=2$ we introduce viscosity solution and prove the convergence of corresponding difference scheme.
 
 \subsection{Notations}
We will make the notations for the one-dimensional and two-dimensional cases parallely. For the sake of simplicity, we will assume that $\Omega=(-1,1)$ in one-dimensional case and $\Omega=(-1,1)\times(-1,1)$ in two-dimensional case in the rest of the paper, keeping in mind that the method works also for more complicated domains.

Let $N\in\mathbb{N}$ be a positive integer, $h=2/N$ and
$$
x_i=-1+ ih,\, y_i=-1+ih,\quad i=0,1,...,N.
$$

We use the notation $u^l_i$ and $u^l_{i,j}$ (or simply $u^l_\alpha$, where $\alpha$ is one- or two-dimensional index) for  the finite difference scheme approximation to  $u_l(x_i)$ and $u_l(x_i,y_j)$, respectively. Concerning the boundary functions $\phi_l,$ we assume they are  extended to be zero everywhere outside the boundary $\partial\Omega,$ for all $l=1,2,\dots, m$. The discrete approximation for these functions will be  $\phi^l_i$ and $\phi^l_{i,j},$ respectively (or simply $\phi^l_\alpha$, where $\alpha$ is one- or two-dimensional index). Note that for the case with two population densities, we additionally will use the following notations: 
\[g_i=\phi^1_i-\phi^2_i=\phi_1(x_i)-\phi_2(x_i)\]
and
\[g_{i,j}=\phi^1_{i,j}-\phi^2_{i,j}=\phi_1(x_i,y_j)-\phi_2(x_i,y_j),\]
in one- and two-dimensional cases, respectively.

In this paper we will use also notations $u=(u_{\alpha})$, $g=(g_{\alpha})$ (not to be confused with functions $u, g$).% Also we will write $(a_{\alpha})\le (b_{\alpha})$ in $\mathcal I$ if $a_{\alpha}\le b_{\alpha}$ for all $\alpha\in\mathcal I$.

Denote
$$
 \mathcal N=\{i:\ 0\leq i\leq N\}\quad\mbox{or}\quad \mathcal N=\{(i,j):\ 0\leq i,j\leq N\},
$$
$$
  \mathcal N^o=\{i:\ 1\leq i\leq N-1\} \quad\mbox{or}\quad  \mathcal N^o=\{(i,j):\ 1\leq i,j\leq N-1\},
$$
in one- and two- dimensional cases, respectively, and
$$
  \partial \mathcal N=\mathcal N \setminus \mathcal N^o.
$$

In one-dimensional case, we consider the following approximation for Laplace operator: for any $i\in \mathcal N^o$,
\[
  \Delta_h v_{i}\equiv L_h v_{i}=\frac{v_{i-1}-2v_{i}+v_{i+1}}{h^2},
\]
and for two-dimensional case we introduce the following 5-point stencil approximation for Laplacian:
\[
 \Delta_h v_{i,j}\equiv L_h v_{i,j}=\frac{v_{i-1,j}+v_{i+1,j}-4v_{i,j}+v_{i,j-1}+v_{i,j+1}}{h^2}
\]
for any $(i,j)\in \mathcal N^o$.

\section{Segregation problem with $m\geq 2$ population densities}

\subsection{The minimizing functional and existence of difference scheme}

In this section, we introduce the discrete counterpart of the spatial segregation problem for the general case. As in the previous section  we assume that $F_l(x,s)$ are convex in the  variable $s,$ for all $l=\overline{1,m},$  and satisfy the properties stated in introduction.

In the rest of the paper the following notation 
\[
\hat{z}^k=z^k-\sum_{j\neq k}z^j,
\]
for elements  $(z^1,z^2,\dots,z^m)$, will play a crucial role.
We focus on the following functional
\begin{multline}\label{disc_energy}
J_h(v^1,v^2,\dots,v^m)=-\frac{1}{2}\sum_{l=1}^m(L_h \hat{v}^l, v^l) +
\sum_{l=1}^m\left(\sum_{\alpha\in\mathcal N}F_l(x_\alpha, v_{\alpha}^l)\right)-\sum_{l=1}^m(L_h \hat{\phi}^l, v^l),
\end{multline}
defined over the set
\begin{multline}\label{disc_min_set}
\mathcal S=\{(v^1,v^2,\dots,v^m)\in (\mathcal H)^{m} :\; v_{\alpha}^p\geq0,\; v_{\alpha}^{p} \cdot v_{\alpha}^{q}=0,\; p\neq q,\; v_{\alpha}^{p}=0 \;\; \text {on} \;\;\partial \mathcal N\},
\end{multline}
where
\[
\mathcal H=\{v=(v_\alpha): v_\alpha\in \mathbb R, \ \alpha\in \mathcal N\}.
\]

Here for $w=(w_\alpha)$ and $v=(v_\alpha)$, $\alpha\in \mathcal N$, the inner product $(\cdot, \cdot)$ is defined by

$$
  (w,v)=\sum_{\alpha\in \mathcal N^0}w_\alpha\cdot v_\alpha.
$$

\begin{theorem}
In view of definitions \eqref{disc_energy} and \eqref{disc_min_set}, the following minimization problem
\begin{equation}\label{disc_min_problem}
   \inf_{\mathcal S}J_h(v^1,v^2,\dots,v^m)\
\end{equation}
has a solution.
\end{theorem}
\begin{proof}
We are going to prove that functional \eqref{disc_energy} is coercive in the set $\mathcal S$. Due to the standard arguments of calculus of variations coercivity and lower semi-continuity will imply the existence of minimizers over the closed set $\mathcal S$. To this end, we observe that 
\[ 
-\frac{1}{2}\sum_{l=1}^m(L_h \hat{v}^l, v^l)=-\frac{1}{2}\sum_{l=1}^m(L_h {v}^l, v^l)+\sum_{1\leq i<j\leq m}(L_h {v}^i, v^j).
\]
For every fixed $i$ and $j$ such that $i\neq j$ we have $v_{\alpha}^i\cdot v_{\alpha}^j=0$. Thus, 
\[
(L_h {v}^i, v^j)=\sum_{\{v_{\alpha}^j>0\}}L_h {v}_{\alpha}^i\cdot v_{\alpha}^j\geq 0
\]
due to the simple fact that $v_{\alpha}^j>0$ implies ${v}_{\alpha}^i=0$, which in turn yields $L_h {v}_{\alpha}^i\geq 0$. Therefore
\[ 
-\frac{1}{2}\sum_{l=1}^m(L_h \hat{v}^l, v^l)\geq-\frac{1}{2}\sum_{l=1}^m(L_h {v}^l, v^l)\geq C\cdot\sum_{l=1}^m(v^l)^2
\]
for some constant $C>0$. Recalling that  $F_l(x,u^l)\geq 0$ we finally obtain that the functional \eqref{disc_energy} is coercive.
\end{proof}

\begin{proposition}\label{onephase}
If an element  $\tilde U=(u^1,u^2,\dots,u^m)\in \mathcal S$  solves the following minimization problem:
\begin{equation*}
   \inf_{\mathcal S}J_h(v^1,v^2,\dots,v^m),\
\end{equation*}
then for every $l\in\{1,2,\dots,m\},$ and $\alpha\in\mathcal N$%  the vector $(w^1,w^2,\dots,w^m)\equiv(u^1+\phi^1,u^2+\phi^2,\dots,u^m+\phi^m)$ solves the following discrete system for every $l=1,2,\dots,m$:
\begin{align}\label{discrete_onphase}
%\begin{cases}
L_h (\hat{u}_\alpha^l+\hat{\phi}_\alpha^l)=f_l(x_\alpha,u_\alpha^l),\;\;\text{whenever}\;\; u_\alpha^l>0.
%\end{cases}
\end{align}

\end{proposition}

\begin{proof}
First of all, it is easy to verify that, if two vectors  $(u^1,u^2,\dots,u^m)$ and $(v^1,v^2,\dots,v^m)$ belong to the set $\mathcal S,$ then
for arbitrary $\varepsilon>0$ we have
\[
((\hat{u}^1+\varepsilon\hat{v}^1)^+,(\hat{u}^2+\varepsilon\hat{v}^2)^+,\dots,(\hat{u}^m+\varepsilon\hat{v}^m)^+)\in \mathcal S,
\]
and
\[
((\hat{u}^1-\varepsilon\hat{v}^1)^+,(\hat{u}^2-\varepsilon\hat{v}^2)^+,\dots,(\hat{u}^m-\varepsilon\hat{v}^m)^+)\in \mathcal S.
\]

Let a vector $(u^1,u^2,\dots,u^m)$  be a minimizer to our discrete minimization problem \eqref{disc_min_problem}. If  $w^i=(\hat{u}^i+\varepsilon\hat{v}^i)^+, \, \varepsilon>0$ we  obtain:
\[
\begin{split}
J_h(w^1,w^2,\dots,w^m)-J_h(u^1,u^2,\dots,u^m)\geq 0,
\end{split}
\]
for every $(v^1,v^2,\dots,v^m)\in \mathcal S$.
Thus,

\begin{multline*}
-\frac{1}{2}\left(\sum_{l=1}^m\left(L_h \hat{w}^l,w^l\right)-\sum_{l=1}^m(L_h \hat{u}^l, u^l)\right)+\\+
\sum_{l=1}^m\left(\sum_{\alpha\in\mathcal N}\left(F_l(x_\alpha, w_{\alpha}^l)-F_l(x_\alpha, u_{\alpha}^l)\right)\right)-\sum_{l=1}^m\left(L_h\hat{\phi}^l, w^l-u^l \right) \geq 0.
\end{multline*}

We choose the vector $(v^1,v^2,\dots,v^m)$ such that for every $l=1,2,\dots,m$ and $\alpha\in\mathcal N$ the following condition holds
\[
\hat{u}_{\alpha}^l\cdot \hat{v}_{\alpha}^l\geq 0.
\]
This implies  the following identity:
$$
 w_{\alpha}^l=(\hat{u}_\alpha^l+\varepsilon\hat{v}_\alpha^l)^+=(\hat{u}_\alpha^l)^++\varepsilon(\hat{v}_\alpha^l)^+={u}_{\alpha}^l +\varepsilon {v}_{\alpha}^l.
$$
Hence,
%\[
\begin{multline}
-\frac{1}{2}\left(\sum_{l=1}^m\sum_{\alpha\in\mathcal N}
{L_h( \hat{u}_{\alpha}^l +\varepsilon\hat{v}_{\alpha}^l)\cdot({u}_{\alpha}^l+\varepsilon{v}_{\alpha}^l)}-\sum_{l=1}^m(L_h \hat{u}^l, u^l)\right)+\\+
\sum_{l=1}^m\left(\sum_{\alpha\in\mathcal N}\left(F_l(x_\alpha, {u}_{\alpha}^l +\varepsilon{v}_{\alpha}^l)-F_l(x_\alpha, u_{\alpha}^l)\right)\right)-\varepsilon\sum_{l=1}^m\left(L_h\hat{\phi}^l, {v}^l\right)=
\end{multline}
\begin{multline}\label{poch}
=-\frac{1}{2}\left(\sum_{l=1}^m\sum_{\alpha\in\mathcal N}L_h \hat{u}_{\alpha}^l\cdot{u}_{\alpha}^l +
\varepsilon^2\sum_{l=1}^m\sum_{\alpha\in\mathcal N}L_h \hat{v}_{\alpha}^l\cdot {v}_{\alpha}^l \right.
+\\+
\left. \varepsilon\sum_{l=1}^m\left(\sum_{\alpha\in\mathcal N}
(L_h \hat{u}_{\alpha}^l\cdot {v}_{\alpha}^l +L_h \hat{v}_{\alpha}^l\cdot {u}_{\alpha}^l+2L_h\hat\phi_{\alpha}^l\cdot {v}_{\alpha}^l) \right)-\sum_{l=1}^m(L_h \hat{u}^l, u^l) \right)+\\+
\sum_{l=1}^m\left(\sum_{\alpha\in\mathcal N}\left(F_l(x_\alpha, {u}_{\alpha}^l +\varepsilon {v}_{\alpha}^l)-F_l(x_\alpha, u_{\alpha}^l)\right)\right)\geq 0.
\end{multline}

%The inequality \eqref{poch} will be reduced to

%\begin{multline}\label{poch_reduced}
%-\frac{1}{2}\left(\varepsilon^2\sum_{l=1}^m\sum_{\alpha\in\mathcal N}L_h {v}_{\alpha}^l\cdot {v}_{\alpha}^l + 2\varepsilon\sum_{l=1}^m\sum_{\alpha\in\mathcal N}L_h ({u}_{\alpha}^l+\phi_{\alpha}^l)\cdot {v}_{\alpha}^l  \right)+\\+
%\sum_{l=1}^m\left(\sum_{\alpha\in\mathcal N}\left(F_l(x_\alpha, {u}_{\alpha}^l +\varepsilon {v}_{\alpha}^l)-F_l(x_\alpha, u_{\alpha}^l)\right)\right)
%\geq 0.
%\end{multline}
 Next, dividing both sides  in \eqref{poch} by $\varepsilon$ and letting $\varepsilon\to 0,$ we arrive at:
 \begin{multline}\label{final_form+}
-\frac{1}{2}\sum_{l=1}^m \sum_{\alpha\in\mathcal N}
\left(L_h\hat{u}_{\alpha}^l\cdot{v}_{\alpha}^l +L_h\hat{v}_{\alpha}^l\cdot{u}_{\alpha}^l+2L_h\hat{\phi}_{\alpha}^l\cdot{v}_{\alpha}^l\right) +
\sum_{l=1}^m\left(\sum_{\alpha\in\mathcal N}f_l(x_\alpha, {u}_{\alpha}^l)\cdot {v}_{\alpha}^l\right)
\geq 0.
\end{multline}

Let $u_{\alpha_0}^{l_0}>0$. We set $(v^1,v^2,\dots,v^m)$ as follows:
\begin{equation}\label{chosen_vector}
\begin{cases}
v_\alpha^l=0,\;\;\mbox{for all}\;\; l\neq l_0,  \alpha\in\mathcal N,\\
v_{\alpha_0}^{l_0}=u_{\alpha_0}^{l_0},\\
v_{\alpha}^{l_0}=0,\;\;\mbox{for all}\;\;  \alpha\in\mathcal N\setminus{\alpha_0}.
\end{cases}
\end{equation}
It is easy to see, that the chosen vector $(v^1,v^2,\dots,v^m)$  satisfies  $\hat{u}_{\alpha}^l\cdot \hat{v}_{\alpha}^l\geq 0$.
Therefore,  the  inequality \eqref{final_form+} holds for this vector, which means it  can be substituted   into   \eqref{final_form+}.
We clearly obtain
\[
-\frac{1}{2}\sum_{l=1}^m \sum_{\alpha\in\mathcal N}
\left(L_h\hat{u}_{\alpha}^l\cdot{v}_{\alpha}^l +L_h\hat{v}_{\alpha}^l\cdot{u}_{\alpha}^l+2L_h\hat{\phi}_{\alpha}^l\cdot{v}_{\alpha}^l\right)=-L_h(\hat{u}_{\alpha_0}^{l_0}+\hat{\phi}_{\alpha_0}^{l_0})\cdot u_{\alpha_0}^{l_0}
\]
and
\[
\sum_{l=1}^m\left(\sum_{\alpha\in\mathcal N}f_l(x_\alpha, {u}_{\alpha}^l)\cdot {v}_{\alpha}^l\right)=f_{l_0}(x_{\alpha_0}, {u}_{\alpha_0}^{l_0})\cdot u_{\alpha_0}^{l_0}.
\]
Hence,
\begin{equation}\label{last_ineq+}
-L_h(\hat{u}_{\alpha_0}^{l_0}+\hat{\phi}_{\alpha_0}^{l_0})+f_{l_0}(x_{\alpha_0}, {u}_{\alpha_0}^{l_0})\geq 0.
\end{equation}

In the same way,  for every $(v^1,v^2,\dots,v^m)\in \mathcal S$ and $\varepsilon>0$ we apparently have:
\[
\begin{split}
J_h((\hat{u}^1-\varepsilon\hat{v}^1)^+,(\hat{u}^2-\varepsilon\hat{v}^2)^+,\dots,(\hat{u}^m-\varepsilon\hat{v}^m)^+)-
J_h(u^1,u^2,\dots,u^m)\geq 0.
\end{split}
\]
In this case we choose  the vector $(v^1,v^2,\dots,v^m)$ such that for every $l=1,2,\dots,m$ and $\alpha\in\mathcal N$ the following condition holds
\[
\hat{v}_{\alpha}^l\cdot(\hat{u}_{\alpha}^l - \hat{v}_{\alpha}^l)\geq 0.
\]
This implies  the following identity:
$$
 (\hat{u}_\alpha^l-\varepsilon\hat{v}_\alpha^l)^+=(\hat{u}_\alpha^l)^+-\varepsilon(\hat{v}_\alpha^l)^+={u}_{\alpha}^l -\varepsilon {v}_{\alpha}^l.
$$

%As in the previous case, for sufficiently small $\varepsilon$   we get  the following identity
%\begin{equation}\label{set_epsilon-}
%  \{\hat{u}_{\alpha}^l-\varepsilon\hat{v}_{\alpha}^l>0\}=\{\hat{u}_{\alpha}^l>0\}\cup\left(\{\hat{u}_{\alpha}^l=0\}\cap\{\hat{v}_{\alpha}^l<0\}\right)\equiv T_{\alpha,l}^-.
%\end{equation}

After proceeding the same steps as above, we obtain
 \begin{multline}\label{final_form-}
 \frac{1}{2}\sum_{l=1}^m \sum_{\alpha\in\mathcal N}
 \left(L_h\hat{u}_{\alpha}^l\cdot{v}_{\alpha}^l +L_h\hat{v}_{\alpha}^l\cdot{u}_{\alpha}^l+2L_h\hat{\phi}_{\alpha}^l\cdot{v}_{\alpha}^l\right) -
 \sum_{l=1}^m\left(\sum_{\alpha\in\mathcal N}f_l(x_\alpha, {u}_{\alpha}^l)\cdot {v}_{\alpha}^l\right)
 \geq 0.
 \end{multline}
%\begin{equation}\label{final_form-}
%\sum_{l=1}^m \sum_{\alpha\in\mathcal N}L_h ({u}_{\alpha}^l+\phi_{\alpha}^l)\cdot {v}_{\alpha}^l -
%\sum_{l=1}^m\left(\sum_{\alpha\in\mathcal N}f_l(x_\alpha, {u}_{\alpha}^l)\cdot {v}_{\alpha}^l\right)
%\geq 0.
%\end{equation}

Here again  we  choose   $(v^1,v^2,\dots,v^m)$ as  in \eqref{chosen_vector}, which apparently satisfies  $\hat{v}_{\alpha}^l\cdot(\hat{u}_{\alpha}^l - \hat{v}_{\alpha}^l)\geq 0$. Therefore, we can substitute the vector  $(v^1,v^2,\dots,v^m)$  into \eqref{final_form-}, which will lead to

\begin{equation}\label{last_ineq-}
L_h (\hat{u}_{\alpha_0}^{l_0}+ \hat{\phi}_{\alpha_0}^{l_0})-f_{l_0}(x_{\alpha_0}, {u}_{\alpha_0}^{l_0})\geq 0.
\end{equation}
Thus, in light of \eqref{last_ineq+} and \eqref{last_ineq-} we obtain
\begin{equation}
L_h (\hat{u}_{\alpha_0}^{l_0}+ \hat{\phi}_{\alpha_0}^{l_0})=f_{l_0}(x_{\alpha_0}, {u}_{\alpha_0}^{l_0}),\;\; \mbox{whenever}\;\; {u}_{\alpha_0}^{l_0}>0.
\end{equation}
\end{proof}

\begin{lemma}\label{aux_lemma}
	Let an element  $(u^1,u^2,\dots,u^m)\in \mathcal S$ solves the minimization problem \eqref{disc_min_problem}, and moreover assume that $u_{\alpha}^{l}=0,$ for some $l=1,2,\dots,m$ and $\alpha\in\mathcal N^o$. Then
$$
L_h(\hat{u}_\alpha^l+\hat{\phi}_\alpha^l)\leq f_{l}(x_{\alpha}, {u}_{\alpha}^{l}).
$$
\end{lemma}
\begin{proof}
	Assume that $u_{\alpha_0}^{l_0}=0,$ for some fixed $\alpha_0\in \mathcal N^o$ and $l_0\in\overline{1,m}$. There are two possibilities: either $\hat{u}_{\alpha_0}^{l_0}=0,$ or $\hat{u}_{\alpha_0}^{l_0}<0$. Let $\hat{u}_{\alpha_0}^{l_0}=0,$ then we take a vector $(v^1,v^2,\dots,v^m)\in \mathcal S$ as follows:
	
	\begin{equation}\label{new_chosen_vector}
		\begin{cases}
			v_\alpha^l=0,\;\;\mbox{for all}\;\; l\neq l_0,  \alpha\in\mathcal N,\\
			v_{\alpha_0}^{l_0}=1,\\
			v_{\alpha}^{l_0}=0,\;\;\mbox{for all}\;\;  \alpha\in\mathcal N\setminus{\alpha_0}.
		\end{cases}
	\end{equation}
	 It is clear that $\hat{u}_{\alpha}^l\cdot \hat{v}_{\alpha}^l\geq 0$ is satisfied  for every $l=1,2,\dots,m$ and $\alpha\in\mathcal N$.  Hence, one can substitute the vector \eqref{new_chosen_vector}  into \eqref{final_form+}. 
	This implies  the inequality \eqref{last_ineq+}, namely
	\[
	L_h(\hat{u}_{\alpha_0}^{l_0}+\hat{\phi}_{\alpha_0}^{l_0})\leq f_{l_0}(x_{\alpha_0}, {u}_{\alpha_0}^{l_0}).
	\]
	Now, if we assume that $\hat{u}_{\alpha_0}^{l_0}<0,$ then there exists some $q_0\neq l_0,$ such  that ${u}_{\alpha_0}^{q_0}>0$. In this case according to Proposition \ref{onephase}, we have 
	\[
	L_h(\hat{u}_{\alpha_0}^{q_0}+\hat{\phi}_{\alpha_0}^{q_0})= 
	f_{q_0}(x_{\alpha_0}, {u}_{\alpha_0}^{q_0})\geq 0.
	\]
	Thus,
	\[
	L_h(\hat{u}_{\alpha_0}^{l_0}+\hat{\phi}_{\alpha_0}^{l_0})=-L_h(\hat{u}_{\alpha_0}^{q_0}+\hat{\phi}_{\alpha_0}^{q_0})-2\sum_{r\neq l_0,q_0}L_h({u}_{\alpha_0}^r+\phi_{\alpha_0}^r)\leq 0\leq f_{l_0}(x_{\alpha_0}, {u}_{\alpha_0}^{l_0}).
	\]
\end{proof}

Next theorem shows the existence of the corresponding finite difference scheme for two or more population densities. Note that for the particular case $F_l(x,u_l)=f_l(x)\cdot u_l(x)$ we obtain the difference scheme for the so-called multi-phase obstacle problem.

\begin{theorem}\label{scheme_exist}
If an element $(u^1,u^2,\dots,u^m)\in\mathcal S$ is a minimizer to \eqref{disc_min_problem}, then the vector $(w^1,w^2,\dots,w^m),$ where  $w_\alpha^l=u_\alpha^l+\phi_\alpha^l,$ solves the  discrete system:
	
	 \begin{equation}\label{scheme_sys}
	 \begin{cases}
	 w_\alpha^{l} =\max \left(\cfrac{-f_{l}(x_\alpha,w_\alpha^{l})h^{2}}{4}+\overline{w}_\alpha^{l} - \sum\limits_{p \neq l}  \overline{w}_\alpha^{p} , \,  0\right),\;\;\alpha\in\mathcal N^o\\
	 w_\alpha^{l}=\phi_\alpha^l,\;\;\alpha\in\mathcal \partial\mathcal N
	 \end{cases}
	 \end{equation}
for every $l=1,2,\dots,m$ and $\alpha\in\mathcal N$. 
Here for a given uniform  mesh on $\Omega\subset \mathbb{R}^2,$   we define     \[
\overline{w}_\alpha^{l}=\frac{1}{4}[w^{l} (x_{i-1},y_{j})+w^{l} (x_{i+1},y_{j})+w^{l}(x_{i},y_{j-1})  +w^{l} (x_{i},y_{j+1})]
\]
as the  average of $w_\alpha^{l}$ for all neighbour points of $\alpha=(x_i,y_j)\in\mathcal N^o$.
\end{theorem}
\begin{proof}

If $(z^1,z^2,\dots,z^m)$ solves the discrete system \eqref{scheme_sys}, then one can easily  see that for every $p,q\in \{1,2,\dots,m\}$ and $\alpha\in \mathcal N^o$ the following property holds:
$$
{z}_\alpha^{p}\cdot {z}_\alpha^{q}=0,\;\; \mbox{whenever}\;\; p\neq q,
$$
which in turn implies that $(z^1-\phi^1,z^2-\phi^2,\dots,z^m-\phi^m)$ belongs to the set $\mathcal S$. Thus the solution to the discrete system \eqref{scheme_sys} itself implies the disjointness property. 

Now, in view of Proposition \ref{onephase} and Lemma \ref{aux_lemma}, we  clearly infer that for every  minimizer $(u^1,u^2,\dots,u^m)\in\mathcal S$ to \eqref{disc_min_problem}, the vector $w_\alpha^l=u_\alpha^l+\phi_\alpha^l,$ satisfies the system \eqref{scheme_sys}. 
\end{proof}
\begin{remark}
We would like to emphasize that the theorem statement holds also for more general uniformly discretized schemes without restriction on the number of stencil (neighbour) points corresponding to a discrete Laplacian. 
\end{remark}

Observe that  we have only used  $f_l(x,s)\geq 0$ to prove the existence of the solution to the discrete system \eqref{scheme_sys}. The convexity assumption on $F_l(x,s)$ will be needed to prove the uniqueness of the discrete scheme which is done in Section \ref{unique_schm}.

\subsection{Uniqueness of difference scheme}\label{unique_schm}
In this section our goal is to show the uniqueness of the difference scheme, which solves the discrete system \eqref{scheme_sys}. To this aim, we need two auxiliary lemmas.

For the sake of convenience  we  denote by $nbr(\alpha)$ the set of all closest neighbour points corresponding to a mesh point $\alpha\in\mathcal N$.  
\begin{lemma}\label{lemma1}
	Let the functions $f_{l}(x,s)$ be nondecreasing with respect to the variable $s$. We take any two elements $(u^{1},u^{2},\dots,u^{m})$ and $(v^{1},v^{2},\dots,v^{m})$ in $\mathcal S$. If  $p_\alpha^l=u_\alpha^l+\phi_\alpha^l,$  and $q_\alpha^l=v_\alpha^l+\phi_\alpha^l,$ are satisfying the discrete system \eqref{scheme_sys}, then the following equation holds:
	\[
	\max_{\mathcal N}(\hat{u}_\alpha^l-\hat{v}_\alpha^l)=\max_{\{u_\alpha^l\leq v_\alpha^l\}}(\hat{u}_\alpha^l-\hat{v}_\alpha^l),
	\]
	for all $l=1,2,\dots,m$.
\end{lemma}

\begin{proof}
	We argue by contradiction. Suppose for some $l_0$ we have
	\begin{equation}\label{init_assmp}
	\hat{u}_{\alpha_0}^{l_0}-\hat{v}_{\alpha_0}^{l_0}=
	\max_{\mathcal N}(\hat{u}_\alpha^{l_0}-\hat{v}_\alpha^{l_0})=\max_{\{u_\alpha^{l_0}> v_\alpha^{l_0}\}}(\hat{u}_\alpha^{l_0}-\hat{v}_\alpha^{l_0})>
	\max_{\{u_\alpha^{l_0}\leq v_\alpha^{l_0}\}}(\hat{u}_\alpha^{l_0}-\hat{v}_\alpha^{l_0}).
	\end{equation}
	Then taking into account the following simple chain of inclusions
	\begin{equation}\label{incl_chain}
		\{u_\alpha^l> v_\alpha^l\}\subset\{\hat{u}_\alpha^l> \hat{v}_\alpha^l\}\subset\{u_\alpha^l\geq v_\alpha^l\},
	\end{equation}
	we obviously see that $ u_{\alpha_0}^{l_0}> v_{\alpha_0}^{l_0}\geq 0 $ implies 
	$\hat{u}_{\alpha_0}^{l_0}>\hat{v}_{\alpha_0}^{l_0}$. On the other hand,  the discrete system \eqref{scheme_sys} gives us  $$L_h(\hat{u}_{\alpha_0}^{l_0}+\hat{\phi}_{\alpha_0}^{l_0})=f_{l_0}(x_{\alpha_0},{u}_{\alpha_0}^{l_0})\;\; \mbox{and}\;\; L_h(\hat{v}_{\alpha_0}^{l_0}+\hat{\phi}_{\alpha_0}^{l_0})\leq f_{l_0}(x_{\alpha_0},{v}_{\alpha_0}^{l_0}).$$ Therefore 
	$$
	L_h(\hat{u}_{\alpha_0}^{l_0}-\hat{v}_{\alpha_0}^{l_0})\geq f_{l_0}(x_{\alpha_0},{u}_{\alpha_0}^{l_0})-f_{l_0}(x_{\alpha_0},{v}_{\alpha_0}^{l_0})\geq 0.
	$$
	Thus,
	$$
	\hat{u}_{\alpha_0}^{l_0}-\hat{v}_{\alpha_0}^{l_0}\leq \frac{1}{4}\sum_{\{\delta\in nbr(\alpha_0)\}}(\hat{u}_{\delta}^{l_0}-\hat{v}_{\delta}^{l_0}),
	$$
	which implies that $\hat{u}_{\alpha_0}^{l_0}-\hat{v}_{\alpha_0}^{l_0}=\hat{u}_{\delta}^{l_0}-\hat{v}_{\delta}^{l_0}>0,$ for all $\delta\in nbr(\alpha_0)$. Due to the chain  \eqref{incl_chain}, we apparently have  ${u}_{\delta}^{l_0}\geq {v}_{\delta}^{l_0}$. According to our assumption \eqref{init_assmp}, the only possibility is ${u}_{\delta}^{l_0}>{v}_{\delta}^{l_0}$  for all $\delta\in nbr(\alpha_0)$. Now we can proceed the previous steps for all $\delta\in nbr(\alpha_0)$ and then for each one we will get corresponding neighbours with the same strict inequality and so on. Continuing this we will finally approach to the boundary $\partial\mathcal N$ where as we know ${u}_{\alpha}^{l_0}={v}_{\alpha}^{l_0}=0$ for all $\alpha\in\partial\mathcal N$. Hence, the strict inequality  fails, which implies that our initial assumption \eqref{init_assmp} is false. Observe that the same arguments can be applied if we interchange the roles of ${u}_{\alpha}^{l}$ and ${v}_{\alpha}^{l}$. Thus, we also have 
	\[
	\max_{\mathcal N}(\hat{v}_\alpha^l-\hat{u}_\alpha^l)=\max_{\{v_\alpha^l\leq u_\alpha^l\}}(\hat{v}_\alpha^l-\hat{u}_\alpha^l),
	\]
	for every $l=1,2,\dots,m$.

	%	It is easy to see that $\max\limits_{\mathcal N}(\hat{v}_\alpha^l-\hat{u}_\alpha^l)=-\min\limits_{\mathcal N}(\hat{u}_\alpha^l-\hat{v}_\alpha^l)$. 
	Particularly, for every fixed $l=1,2.\dots,m$ and $\alpha\in\mathcal N$ we have
	\begin{multline}\label{double_ineq}
	-\max_{\{v_\alpha^l\leq u_\alpha^l\}}(\hat{v}_\alpha^l-\hat{u}_\alpha^l)=
	-\max\limits_{\mathcal N}(\hat{v}_\alpha^l-\hat{u}_\alpha^l)\leq \hat{u}_\alpha^l-\hat{v}_\alpha^l\leq \\ \le 
	\max\limits_{\mathcal N}(\hat{u}_\alpha^l-\hat{v}_\alpha^l)=\max_{\{u_\alpha^l\leq v_\alpha^l\}}(\hat{u}_\alpha^l-\hat{v}_\alpha^l).
	\end{multline}
\end{proof}
Thanks to Lemma \ref{lemma1} in the sequel we will use the following notations:
$$
M:=\max_l\;\left(\max_{\mathcal N}\left(\hat{u}_\alpha^l-\hat{v}_\alpha^l\right)\right)=\max_l\;\left(\max_{\{u_\alpha^l\leq v_\alpha^l\}}\left(\hat{u}_\alpha^l-\hat{v}_\alpha^l\right)\right),
$$
and
$$
R:=\max_l\;\left(\max_{\mathcal N}\left(\hat{v}_\alpha^l-\hat{u}_\alpha^l\right)\right)=\max_l\;\left(\max_{\{v_\alpha^l\leq u_\alpha^l\}}\left(\hat{v}_\alpha^l-\hat{u}_\alpha^l\right)\right).
$$

\begin{lemma}\label{lemma2}
	Let the functions $f_{l}(x,s)$ be nondecreasing with respect to the variable $s$. We also assume that for   two elements $(u^{1},u^{2},\dots,u^{m})$ and $(v^{1},v^{2},\dots,v^{m})$ in $\mathcal S$ we have  $p_\alpha^l=u_\alpha^l+\phi_\alpha^l,$  and $q_\alpha^l=v_\alpha^l+\phi_\alpha^l,$ are solving the discrete system \eqref{scheme_sys}. For these two elements we set $M$ and $R$ as defined above. If $M>0$ and it is attained for some $l_0$, then $M=R>0$  and there exists some $t_0\neq l_0,$ and $\delta_0\in\mathcal{N},$ such that
$$
0<M=\max_{\{u_\alpha^{l_0}\leq v_\alpha^{l_0}\}}(\hat{u}_\alpha^{l_0}-\hat{v}_\alpha^{l_0})=\max_{\{u_\alpha^{l_0}= v_\alpha^{l_0}=0\}}(\hat{u}_\alpha^{l_0}-\hat{v}_\alpha^{l_0})={v}_{\delta_0}^{t_0}-{u}_{\delta_0}^{t_0}.
$$
	
\end{lemma}
\begin{proof}
	It is easy to verify that $(\hat{u}_\alpha^{l_0}-\hat{v}_\alpha^{l_0})$ might be positive only on the set $\{u_\alpha^{l_0}= v_\alpha^{l_0}=0\}$ (for the other cases $(\hat{u}_\alpha^{l_0}-\hat{v}_\alpha^{l_0})\leq 0)$. Hence, 
		$$
		M=\max_{\{u_\alpha^{l_0}= v_\alpha^{l_0}=0\}}(\hat{u}_\alpha^{l_0}-\hat{v}_\alpha^{l_0}).
		$$
	Using the latter equality, one can prove that $M>0$ implies $R>0$.  Indeed, if we assume that $R\leq 0,$ then according to definition of $R$ we will get that $\hat{v}_\alpha^{l}\leq\hat{u}_\alpha^{l}$ for all $l=\overline{1,m}$ and $\alpha\in\mathcal N$. This obviously yields ${v}_\alpha^{l}\leq {u}_\alpha^{l},$ for all $l=\overline{1,m}$ and $\alpha\in\mathcal N$.  Thus, 
	$$
	0<M=\max_{\{u_\alpha^{l_0}= v_\alpha^{l_0}=0\}}(\hat{u}_\alpha^{l_0}-\hat{v}_\alpha^{l_0})=
	\max_{\{u_\alpha^{l_0}= v_\alpha^{l_0}=0\}}\left(
	\sum_{l\neq l_0}(v_{\alpha}^{l}-u_{\alpha}^{l})\right)\leq 0.
	$$
	This is a contradiction, and therefore $R>0$. It is apparent that in a similar way one can prove the converse statement as well. Thus, we clearly see that at the same time either both $R$ and $M$ are non-positive, or they are positive.

Concerning the equality $M=R,$ it is easy to see that if the maximum $M>0$ is attained at the mesh point $\delta_0\in\mathcal{N},$  then the following holds:
	\begin{equation*}
	0<M=\max_{\{u_\alpha^{l_0}= v_\alpha^{l_0}=0\}}(\hat{u}_\alpha^{l_0}-\hat{v}_\alpha^{l_0})=
	\hat{u}_{\delta_0}^{l_0}-\hat{v}_{\delta_0}^{l_0}=\sum_{l\neq l_0}(v^l_{\delta_0}-u^l_{\delta_0}).
	\end{equation*}
Since  $\sum_{l\neq l_0}v^l_{\delta_0}$ is positive,  then there exists  $t_0\neq l_0$ such that ${v}_{\delta_0}^{t_0}>0$. Thus, we conclude
	\begin{equation}\label{M=R}
		0<M={v}_{\delta_0}^{t_0}-\sum_{l\neq l_0}{u}_{\delta_0}^{l}\leq\hat{v}_{\delta_0}^{t_0}-u_{\delta_0}^{t_0}+\sum_{l\neq t_0}u^l_{\delta_0}= \hat{v}_{\delta_0}^{t_0}-\hat{u}_{\delta_0}^{t_0}\leq R.
	\end{equation}
	In the same way we will obtain that $R\leq M,$ and therefore $M=R$. On the other hand, since $M=R,$ then the following obvious inequality holds
	$$
	M={v}_{\delta_0}^{t_0}-\sum_{l\neq l_0}{u}_{\delta_0}^{l}\geq \hat{v}_{\delta_0}^{t_0}-\hat{u}_{\delta_0}^{t_0}.
	$$ 
	This leads to $2\sum\limits_{l\neq t_0}{u}_{\delta_0}^{l}\leq 0,$ and therefore
	${u}_{\delta_0}^{l}=0,$ for all $l\neq t_0$. Hence,
	$$M={v}_{\delta_0}^{t_0}-\sum_{l\neq l_0}{u}_{\delta_0}^{l}={v}_{\delta_0}^{t_0}-{u}_{\delta_0}^{t_0}.
	$$
\end{proof}

Now, we are in a position to prove the main result of this section.
\begin{theorem}\label{scheme_unique}
	Let the functions $f_{l}(x,s)$ are nondecreasing with respect to the variable $s$. Then
	there exists  a unique vector  $(u^{1},u^{2},\dots,u^{m})\in\mathcal S,$ such that
	$w_\alpha^{l}=u_\alpha^{l}+\phi_\alpha^{l}$  satisfies the discrete system \eqref{scheme_sys}.
\end{theorem}

\begin{proof}

	Let two elements $(u^{1},u^{2},\dots,u^{m})$ and $(v^{1},v^{2},\dots,v^{m})$ in $\mathcal S$ fulfil the theorem statement, i.e.  $(u^{1}+\phi^1,u^{2}+\phi^2,\dots,u^{m}+\phi^m)$ and $(v^{1}+\phi^1,v^{2}+\phi^2,\dots,v^{m}+\phi^m)$ be solving the discrete system \eqref{scheme_sys}. For these vectors we set the definition of $M$ and $R$. Then, we consider two cases $M\leq 0$ and $M>0$. If we assume that $M\leq 0,$ then according to Lemma \ref{lemma2}, we get  $R\leq 0$. But if  $M$ and $R$ are non-positive, then the uniqueness follows. Indeed,  due to  \eqref{double_ineq} we have the following obvious inequalities
	$$
	0\leq -R\leq \hat{u}_\alpha^l-\hat{v}_\alpha^l\leq M\leq 0.
	$$
	This provides  for every $l=\overline{1,m}$ and $\alpha\in\mathcal N$ we have $\hat{u}_\alpha^l=\hat{v}_\alpha^l,$ which in turn implies  $${u}_\alpha^l={v}_\alpha^l.$$
	Now suppose  $M>0$. Our aim is to prove that this case leads to a contradiction. Let the value $M$ is attained for some $l_0\in\overline{1,m},$ then
	due to Lemma \ref{lemma2} there exist $\delta_0\in \mathcal N$ and $t_0\neq l_0$ such that:
	$$0<M=R=\max_{\{u_\alpha^{l_0}\leq v_\alpha^{l_0}\}}(\hat{u}_\alpha^{l_0}-\hat{v}_\alpha^{l_0})=\max_{\{u_\alpha^{l_0}= v_\alpha^{l_0}=0\}}(\hat{u}_\alpha^{l_0}-\hat{v}_\alpha^{l_0})={v}_{\delta_0}^{t_0}-{u}_{\delta_0}^{t_0}.$$
	Using the fact that $f_{l}(x,s)$ are nondecreasing with respect to the variable $s$, we clearly obtain 
	$$
	L_h(\hat{v}_{\delta_0}^{t_0}-\hat{u}_{\delta_0}^{t_0})\geq 0.
	$$
	Thus, 
	$$
	\hat{v}_{\delta_0}^{t_0}-\hat{u}_{\delta_0}^{t_0}\leq \frac{1}{4}\sum_{\{\gamma\in nbr(\delta_0)\}}(\hat{v}_{\gamma}^{t_0}-\hat{u}_{\gamma}^{t_0}),
	$$
	which implies  $M=\hat{v}_{\delta_0}^{t_0}-\hat{u}_{\delta_0}^{t_0}=\hat{v}_{\gamma}^{t_0}-\hat{u}_{\gamma}^{t_0}>0$ for all $\gamma\in nbr(\delta_0)$.  The chain \eqref{incl_chain} provides that  for all $\gamma\in nbr(\delta_0),$ we have ${v}_{\gamma}^{t_0}\geq{u}_{\gamma}^{t_0}$. For the neighbour mesh points $\gamma$ we proceed as follows: If  ${v}_{\gamma}^{t_0}>{u}_{\gamma}^{t_0}$ for some $\gamma_0\in nbr(\delta_0),$ then obviously $L_h(\hat{v}_{\gamma_0}^{t_0}-\hat{u}_{\gamma_0}^{t_0})\geq 0$. 
	This, as we saw a few lines above, leads to  $M=\hat{v}_{\gamma_0}^{t_0}-\hat{u}_{\gamma_0}^{t_0}=\hat{v}_{\theta}^{t_0}-\hat{u}_{\theta}^{t_0}>0$ for all $\theta\in nbr(\gamma_0)$. 
	
 If  ${v}_{\gamma}^{t_0}={u}_{\gamma}^{t_0}$ for some $\gamma_0\in nbr(\delta_0),$ then due to $\hat{v}_{\gamma_0}^{t_0}-\hat{u}_{\gamma_0}^{t_0}=M>0$ one has ${v}_{\gamma_0}^{t_0}={u}_{\gamma_0}^{t_0}=0$. Hence there exists some $\lambda_0\neq t_0,$ such that
	$$
	M=\hat{v}_{\gamma_0}^{t_0}-\hat{u}_{\gamma_0}^{t_0}=
	\sum_{l\neq t_0}\left({u}_{\gamma_0}^{l}-{v}_{\gamma_0}^{l} \right)={u}_{\gamma_0}^{\lambda_0}-\sum_{l\neq t_0}{v}_{\gamma_0}^{l}.
	$$
	As before, we can write the following inequality
	$$
	M={u}_{\gamma_0}^{\lambda_0}-\sum_{l\neq t_0}{v}_{\gamma_0}^{l}\geq \hat{u}_{\gamma_0}^{\lambda_0}-\hat{v}_{\gamma_0}^{\lambda_0},
	$$
	which in turn gives $2\sum\limits_{l\neq \lambda_0}{v}_{\gamma_0}^{l}\leq 0,$ and therefore
	${v}_{\gamma_0}^{l}=0$ for all $l\neq \lambda_0$. Hence
	$$
	M={u}_{\gamma_0}^{\lambda_0}-\sum_{l\neq t_0}{v}_{\gamma_0}^{l}={u}_{\gamma_0}^{\lambda_0}-{v}_{\gamma_0}^{\lambda_0}.
	$$
	This suggests that we can apply the same approach as above and using the fact that $L_h(\hat{u}_{\gamma_0}^{\lambda_0}-\hat{v}_{\gamma_0}^{\lambda_0})\geq 0,$ one gets that $M=\hat{u}_{\gamma_0}^{\lambda_0}-\hat{v}_{\gamma_0}^{\lambda_0}=\hat{u}_{\theta}^{\lambda_0}-\hat{v}_{\theta}^{\lambda_0}>0,$ for all $\theta\in nbr(\gamma_0)$.
	Thus, all the time continuing this process for the neighbour points, we observe that for every mesh point $\gamma$ there always exists some $l_\gamma\in\overline{1,m}$ such that 
	 $$\mbox{either}\;\;
	\hat{u}_{\gamma}^{l_\gamma}-\hat{v}_{\gamma}^{l_\gamma}=M, \;\;\mbox{or}\;\; 
	\hat{u}_{\gamma}^{l_\gamma}-\hat{v}_{\gamma}^{l_\gamma}=-M.
	$$
	On the other hand, it is clear that sooner or later, we will reach the boundary $\partial\mathcal N,$ and this will give a contradiction, because for every $\gamma\in\partial\mathcal N,$ and $l=\overline{1,m}$ one has $\hat{u}_{\gamma}^{l}-\hat{v}_{\gamma}^{l}=\hat{v}_{\gamma}^{l}-\hat{u}_{\gamma}^{l}=0$.  
\end{proof}

\begin{corollary}
Note that due to  Theorems \ref{scheme_exist} and \ref{scheme_unique},  the solution to the  minimization problem \eqref{disc_min_problem}  is unique, provided all functions $f_{l}(x,s)$ are nondecreasing (all functions $F_{l}(x,s)$ are convex) with respect to the variable $s$.  
\end{corollary}

\begin{remark}
 It is noteworthy that the system \eqref{scheme_sys}, for the case $m=2,$ where $f_{l}(x,s)=f_{l}(x),\; l=1,2,$  has already been suggested in \cite{Farid} to approximate the solution to a two-phase obstacle (membrane) problem. The convergence of the iterative algorithm corresponding to this particular case  was proved in \cite{avetikconv}.
\end{remark}

%
%In view of Section  \ref{twopop_section} we clearly see  two equivalent ways to define a difference method for the segregation problem with two competing densities. The  first  approach is to consider the solution to a discrete Min-Max equation \eqref{nonlinprob}, and the second approach is the solution to a discrete system \eqref{scheme_sys_two}. 
%	
%}

%$\Omega_h=\{\alpha\cdot h: \alpha\in\mathcal N^o\}$. Let also $\partial \Omega_h=\{\alpha\cdot h: \alpha\in\partial \mathcal N\}$.

\subsection{Numerical algorithm and its properties}

Once we know that the solution to a discrete system \eqref{scheme_sys} is unique, we can start to elaborate a numerical algorithm corresponding to that difference scheme. To this end, we suggest the  generalized version  of the  algorithm developed in \cite{MR2563520,Mywork}. The  iterative method for the case of arbitrary $m$  competing densities is  defined as follows:
 \begin{itemize}
 \item \textbf{ Initialization:}
$\text{For}\quad    l=1,\cdots, m, \;\; \text{set}$
\begin{equation*}
(u_{\alpha}^{l})^{(0)}=
\left \{
\begin{array}{ll}
  0  &    \alpha \in \mathcal N^{\circ},  \\
  \phi_{\alpha}^l  &   \alpha\in \partial \mathcal N.
\end{array}
\right.
\end{equation*}
\item  \textbf{{Step $k+1 $, $k\geq 0$:}}
For $l=1,\cdots, m,$ we iterate  for all  interior points
  \begin{equation}\label{algorithm}
 (u_\alpha^{l})^{(k+1)} =\max \left(\cfrac{-f_{l}(x_\alpha,(u_\alpha^{l})^{(k)})h^{2}}{4}+(\overline{u}_\alpha^{l})^{(k)} - \sum_{p \neq l}  (\overline{u}_\alpha^{p})^{(k)} , \,  0\right).
  \end{equation}
  \end{itemize}

\begin{lemma}\label{disjoint}
The iterative algorithm \eqref{algorithm} satisfies
\[(u_\alpha^l)^{(k)}\cdot (u_\alpha^q)^{(k)}=0,\]
for all $k\in\mathbb{N},\alpha\in\mathcal N^o,$ and $q,l\in\{1,2,\dots,m\},where\;\; q\neq l$.
\end{lemma}
\begin{proof}
The proof repeats much the same lines as in \cite[Lemma $2.7$]{Mywork}.
\end{proof}

\begin{lemma}\label{lemma_2}
Numerical algorithm \eqref{algorithm} is stable and consistent.
\end{lemma}
\begin{proof}
Here, we will prove the stability of the method. The proof of the consistency  is straightforward.
Due to $f_l(x,u^l)\geq 0,$ for every $l=1,2,\dots,m,$ we can write the following inequality:
\[
(u_\alpha^{l})^{(k+1)} =\max \left(\cfrac{-f_{l}(x_\alpha,(u_\alpha^{l})^{(k)})h^{2}}{4}+(\overline{u}_\alpha^{l})^{(k)} - \sum_{p \neq l}  (\overline{u}_\alpha^{p})^{(k)} , \,  0\right)
\leq (\overline{u}_\alpha^{l})^{(k)}.
\]

Therefore
\[
(u_\alpha^{l})^{(k+1)}\leq (\overline{u}_\alpha^{l})^{(k)},
\]
which is the same as 
\[
 L_h (u_\alpha^{l})^{(k+1)} \geq 0,
\]
where $L_h$ is the  discrete Laplace operator. After applying the discrete maximum principle we obtain
\[0\leq (u_\alpha^{l})^{(k+1)}\leq\max_{\alpha}\phi_\alpha^{l}.
\]
Hence, $(u_\alpha^{l})^{(k)}$  is uniformly bounded  for every $k\in\mathbb{N}$. 
\end{proof}

\section{The special case $m=2$ revisited }\label{twopop_section}

In this section we revisit the special case $m=2$. Our intention is to obtain the convergence of the difference scheme for this particular case  by introducing the notion of viscosity solution.  This approach cannot be applied for three or more population densities. 

\subsection{Generalized two-phase obstacle problem}

 Let the functions $F_i(x,s),\;i=1,2,$ be convex in the  variable $s,$ and satisfy the properties given in introduction.  Then, employing the  same analysis as in \cite[Section $2.1$]{Mywork}, one can treat the problem \eqref{m=2}  as a  minimization problem subject to the closed and convex set. The rewritten minimization problem will be the following:
\begin{equation}\label{eq-1}
\mbox{Minimize:}\;\;\mathcal J(w)=\int_\Omega \left( \frac{1}{2} |\nabla w|^2 + F_1(x,w^+) + F_2(x,-w^-)\right) dx,
\end{equation}
over the set  $\mathbb K=\{ w\in H^1(\Omega):\;\; w-(\phi_1-\phi_2) \in H^1_0(\Omega)\},$ where  $w^+=\max(w,0)$, $w^-=\min(w,0)$.
If we denote $v_1=w^+$ and $v_2=-w^-,$ then $(v_1,v_2)\in\mathcal{S},$  and it is a unique minimizer to the problem \eqref{m=2}.

The corresponding  Euler-Lagrange equation for the minimization problem  will be
\begin{equation}\label{2D-Problem-DifForm}
\left\{
\begin{array}{ll}
\Delta w = f_1(x,w)\cdot\chi_{\{w>0\}}-f_2(x,-w)\cdot\chi_{\{w<0\}}, & x\in\Omega,\\
w=\phi_1-\phi_2, & x\in\partial \Omega,
\end{array}
\right.
\end{equation}
where $\chi_A$ stands for the characteristic function of the set $A$. Inspired by the setting of the two-phase obstacle problem (see \cite{PSU2012}) we will call the problem \eqref{2D-Problem-DifForm} the \emph{generalized two-phase obstacle problem}. Nowadays, the theory of the two-phase obstacle-like problems (elliptic and parabolic versions) is well-established  and for a reference we again address to the book \cite{PSU2012} and references therein. The interested reader is also referred to the recent works \cite{faridperturb,repin2015}.   For the numerical treatment of the same problems we refer to the works \cite{arakelyan2015finite,ABP2014,MR2961456,Farid,Osher}.

In \cite{ABP2014,MR2961456} the authors introduced the so-called  Min-Max formulation for the usual two-phase obstacle problem, which is very useful to define the notion of viscosity solution. Moreover,  it turns out that the  viscosity solution is equivalent to the weak solution to the two-phase obstacle problem.  Our aim is to use the same approach to the generalized counterpart.
To this end, we need to make some notations.

Let $\Omega$ be an open subset of $\mathbb{R}^n$, and for a twice differentiable function $u:\Omega\to\mathbb R$ let $Du$ and $D^2 u$ denote the gradient and Hessian matrix of $u$, respectively. Also let the function $G(x, r, p, X)$ be a continuous real-valued function defined on $\Omega\times\mathbb{R}\times\mathbb{R}^n\times S^n$, with $S^n$ being the space of real symmetric $n\times n$ matrices.

In light of the Min-Max form defined in \cite{ABP2014} we introduce the following generalized Min-Max variational equation:
\begin{equation}\label{2D-Problem-Viscosity}
\left\{
\begin{array}{ll}
\min\left(-\Delta w + f_1(x,w) ,\max(-\Delta w-f_2(x,-w),w)\right)=0,& \mbox{in}\ \Omega\\
w=\phi_1-\phi_2\equiv g,& \mbox{on}\ \partial\Omega.
\end{array}
\right.
\end{equation}

According to the above notations, we introduce a function $G:\Omega\times\mathbb{R}\times\mathbb{R}^n\times S^n\to\mathbb R$ by
\begin{equation}\label{Matrix_form}
G(x, r, p,X)=\min(-trace(X) + f_1(x,r) ,\max(-trace(X)- f_2(x,-r),r)),
\end{equation}
then the equation in \eqref{2D-Problem-Viscosity} can be rewritten as
\begin{equation}\label{2D-F}
G(x, w, Dw, D^2 w)=0\quad \mbox{in}\ \ \Omega.
\end{equation}

Below we recall the definition of degenerate ellipticity and prove that the equation \eqref{2D-Problem-Viscosity} is degenerate elliptic.

\begin{definition}
We call the equation \eqref{2D-F}  \textbf{degenerate elliptic} if
$$
 F(x, r, p,X)\leq F(x, s, p, Y )  \quad\mbox{whenever }\quad  r\leq s  \quad\mbox{and }\quad Y \leq X,
$$
where $Y \leq X$ means that $X - Y$ is a nonnegative definite symmetric matrix.
\end{definition}
\begin{lemma}\label{2D-lem-degelliptic}
Equation \eqref{Matrix_form} is degenerate elliptic.
\end{lemma}
\begin{proof} Let $X,Y\in S^n$ and $r,s\in\mathbb R$ satisfy $Y \leq X,$ and $r\leq s$. Then using the fact that $F_i(x,t)$ is convex in $t$ for all $x\in\Omega,$  we have $F''_{i}(x,t)=f'_{i}(x,t)\geq 0,$ where the derivatives are taken  with respect to $t$.
Thus,
\begin{equation*}
-trace(X)+f_1(x,r) \leq-trace(Y)+f_1(x,s), 
\end{equation*}
and
\begin{equation*}
\max(-trace(X)-f_2(x,-r),r)\le \max(-trace(Y)-f_2(x,-s),s).
\end{equation*}
Therefore
\begin {align*}
G(x, r, p,X) &=\min(-trace(X) + f_1(x,r) ,\max(-trace(X)-f_2(x,-r),r))\\
             &\leq\min(-trace(Y) +f_1(x,s) ,\max(-trace(Y)-f_2(x,-s),s))\\&=G(x, s, p,Y).
\end{align*}
\end{proof}

Now, we are ready to define viscosity solutions for the generalized two-phase obstacle problem.  For general background about the theory of viscosity solutions the reader is referred to \cite{MR1118699}, \cite{MR1351007} and references therein.
\begin{definition}\label{visc_definition}
A bounded uniformly continuous function $w:\Omega\to\mathbb R$ is called  a viscosity subsolution (resp. supersolution) for \eqref{2D-Problem-Viscosity}, if for each $\varphi\in C^2(\Omega)$ and local maximum point of $w-\varphi$ (respectively minimum) at $x_0\in \Omega,$  we have
\[
\min\left(-\Delta \varphi(x_0) + f_1(x_0,\varphi(x_0)) ,\max(-\Delta \varphi(x_0)-f_2(x_0,-\varphi(x_0)),\varphi(x_0))\right)\leq 0.
\]
(\text{respectively}
\[
\min\left(-\Delta \varphi(x_0) + f_1(x_0,\varphi(x_0)) ,\max(-\Delta \varphi(x_0)-f_2(x_0,-\varphi(x_0)),\varphi(x_0))\right)\geq 0.)
\]

\end{definition}

The function $w:\Omega\to\mathbb R$ is called  a viscosity solution to \eqref{2D-Problem-Viscosity}, if simultaneously it is  a viscosity subsolution and supersolution  for \eqref{2D-Problem-Viscosity}.

%\begin{definition}
%$w:\Omega\to\mathbb R$ is called a \textbf{viscosity subsolution} of \eqref{2D-Problem-Viscosity}, if it is upper semicontinuous and for each $\varphi\in C^2(\Omega)$ and local maximum point $x_0\in \Omega$ of $w-\varphi$ we have
%\[
%\min\left(-\Delta \varphi(x_0) + f_1(x_0,w^+(x_0)) ,\max(-\Delta \varphi(x_0)-f_2(x_0,-w^-(x_0)),w(x_0))\right)\leq 0.
%\]
%\end{definition}

%\begin{definition}
%$w:\Omega\to\mathbb R$ is called a \textbf{viscosity supersolution} of \eqref{2D-Problem-Viscosity}, if it is lower semicontinuous and for each $\varphi\in C^2(\Omega)$ and local minimum point $x_0\in \Omega$ of $w-\varphi$ we have
%\[
%\min\left(-\Delta \varphi(x_0) + f_1(x_0,w^+(x_0)) ,\max(-\Delta \varphi(x_0)-f_2(x_0,-w^-(x_0)),w(x_0))\right)\geq 0.
%\]
%\end{definition}
%\begin{definition}
%$w:\Omega\to\mathbb R$ is called a \textbf{viscosity solution} of \eqref{2D-Problem-Viscosity}, if it is both a viscosity subsolution and supersolution (and hence continuous) for \eqref{2D-Problem-Viscosity}.
%\end{definition}

\subsection{Convergence of finite difference scheme}

Here, according to \cite{ABP2014} and generalized Min-Max variational equation \eqref{2D-Problem-Viscosity}, we define appropriate finite difference scheme as follows:

\begin{equation}\label{nonlinprob}
\begin{cases}
\min(-L_h u_{\alpha}+f_1(x_\alpha, u_{\alpha}) ,\max(-L_h u_{\alpha}-f_2(x_\alpha, -u_{\alpha}) ,  u_{\alpha}))=0, \\[6pt]
u_{\alpha}=g_{\alpha}, \quad\alpha \in \partial \mathcal N.
\end{cases}
\end{equation}
The existence and uniqueness of the solution to the system \eqref{nonlinprob} implicitly follows from the previous section of the paper.   It also can  be shown  directly by  considering the minimization of the following functional (see \cite{ABP2014}):
\begin{equation}\label{disc_eq-1}
J_h(v)=-\frac{1}{2} \Big(L_h v,v\Big)+\sum_{\alpha\in \mathcal N}F_1(x_\alpha, v_{\alpha}^+) +
\sum_{\alpha\in \mathcal N}F_2(x_\alpha, -v_{\alpha}^-) - \Big(L_h g,v\Big),
\end{equation}
subject to the finite dimensional space
$$
\mathcal K=\{v\in\mathcal H: \ v_{\alpha}=0, \ \alpha\in\partial \mathcal N\},\quad {\rm where}\quad
\mathcal H=\{v=(v_\alpha): v_\alpha\in \mathbb R, \ \alpha\in \mathcal N\}.
$$
Here $w=(w_\alpha)$ and $v=(v_\alpha)$, $\alpha\in \mathcal N$. The inner product $(\cdot, \cdot)$ is defined by
$$
(w,v)=\sum_{\alpha\in \mathcal N^0}w_\alpha\cdot v_\alpha.
$$

Our next step is to apply Barles-Souganidis theorem (see \cite{MR1115933}) for viscosity solutions to obtain the convergence of a difference scheme, which solves the system \eqref{nonlinprob}. We will prove the result  for the general uniform structured discretization. Our goal is to show that there is no restriction on the number of stencil points for a discrete Laplacian, as long as  the finite difference schemes with uniform grid are concerned.   To this aim,  we define a uniform  structured grid on the domain $\Omega$ as a directed graph consisting of a set of points $x_i \in \Omega$, $i = 1,\dots,N,$   each endowed with a number of neighbours $K$. A grid function is a real valued function defined on the grid, with values $u_i = u(x_i)$. The typical examples of such grid are $3$-point and $5$-point stencil discretization for the spaces of one dimension and two dimension, respectively. Here we define  degenerate elliptic schemes introduced by Oberman (see \cite {MR2218974}).

%where $\mathcal{L}_h^1u_\alpha=-L_h u_{\alpha}+f_1(x_\alpha, u_{\alpha})$, and  $\mathcal{L}_h^2u_\alpha=-L_h u_{\alpha}-f_2(x_\alpha, -u_{\alpha})$.

%\bblue{\begin{remark}
	
%It is easy to see that for the case $m=2$ the functional \eqref{disc_energy} studied in previous section  will be reduced to 
%\begin{multline*}
%J_h(v^1,v^2)=-\frac{1}{2} \Big(L_h (v^1-v^2),v^1-v^2\Big)+\\+ \sum_{\alpha\in \mathcal N}F_1(x_\alpha, v_{\alpha}^1) +
%\sum_{\alpha\in \mathcal N}F_2(x_\alpha, v_{\alpha}^2) - \Big(L_h (\phi^1-\phi^2),v^1-v^2\Big).
%\end{multline*}
%If we set $u=v^1-v^2$ and $g=\phi^1-\phi^2,$ then we will derive to the two-phase functional \eqref{disc_eq-1} .
		
%\end{remark}
%}

A function $F^h:\mathbb{R}^N\rightarrow\mathbb{R}^N,$ which is regarded as a map from grid functions to grid functions, is a \emph{finite difference scheme} if
\[F^h[u]^i = F^i[u_i, u_i-u_{i_1},\dots,u_i-u_{i_{K}}] \quad (i = 1,\dots,N),\]
where $\{i_1,i_2,\dots,i_{K}\}$ are the neighbour points of a grid point $i$. Denote
\[
F^i[u]\equiv F^i[u_i,u_i-u_{i_j}|_{j=\overline{1,K}}]\equiv F^i[u_i,u_i-u_j],\;\;  i = 1, . . . ,N ,
\]
where $u_j$ is shorthand for the list of neighbours $u_{i_j}|_{j=\overline{1,K}}$.

\begin{definition}\label{def_deg_elliptic_scheme}
 The scheme $F$ is \emph{degenerate elliptic} if each component $F^i$ is nondecreasing
in each variable, i.e. each component of the scheme $F^i$ is a nondecreasing
function of $u_i$ and the differences $u_i - u_{i_j}$ for $j = 1, . . . ,K$.
\end{definition}

Since the grid is uniformly structured, we denote  $h>0$ as the size of the mesh. Then, for the nonlinear system \eqref{nonlinprob} we have
\begin{equation}\label{numscheme}
F^i[u_i,u_i-u_j]= \min(-L_h u_i + f_1(x_i, u_i)\, ,\, \max(-L_h u_i- f_2(x_i, -u_i)\, ,\,  u_i)),
\end{equation}
where
\begin{equation}
L_h u_i=\sum_{j=1}^{K}\frac{1}{h^2}(u_{i_j}-u_{i}),\;\;  i = 1, . . . ,N .
\end{equation}
Since the functions $f_i(x,s)$ are monotone nondecreasing  with respect to $s,$ we clearly see that $F^i[u_i,u_i-u_j]$ is non-decreasing with respect to $u_i$ and $u_i-u_j$ as well. Therefore, the finite difference scheme \eqref{nonlinprob} is a \emph{degenerate elliptic scheme}. We know that the degenerate elliptic schemes are \emph{monotone} and \emph{stable} (see \cite {MR2218974}).
The consistency of the system \eqref{nonlinprob} is obvious. Thus, we show that the nonlinear system  \eqref{nonlinprob} fulfils  all the necessary properties for the Barles-Souganidis framework, namely it is stable, monotone, and consistent, and therefore, due  to Barles - Souganidis theorem, the solution to the discrete nonlinear system \eqref{nonlinprob} converges locally uniformly to the unique viscosity solution to the generalized two-phase obstacle problem \eqref{2D-Problem-Viscosity}.

We summarize the above result in the following theorem.
\begin{theorem}[Convergence]
The finite difference scheme given by system \eqref{numscheme} converges uniformly on compacts subsets of  $\Omega$ to the unique viscosity solution to generalized two phase obstacle equation \eqref{2D-Problem-Viscosity}.
\end{theorem}

\section{Numerical examples}

In this section we present numerical simulations for  two and more competing densities with different internal dynamics  $F_i(x,u_i(x))$.  We consider the following minimization problem:
\begin{equation}\label{last}
\text{Minimize} \int_{\Omega}  \sum_{i=1}^{m} \left( \frac{1}{2}| \nabla u_{i}|^{2}+F_i(x,u_i(x))  \right) dx,
\end{equation}
over the set
$$
S={\{(u_1,\dots,u_{m})\in (H^{1}   (\Omega))^{m} :u_{i}\geq0, u_{i} \cdot u_{j}=0, u_{i}=\phi_{i} \quad \text {on} \quad \partial  \Omega}\}.
$$

\begin{example}
  We take $\Omega=[-1,1]\times[-1,1]$, and the  internal dynamics such that 
  \[f_1(x,y,u_1)=2 (x^2 + y^2) |u_1|,\;\; \mbox{and}\;\; f_2(x,y,u_2)=10 (x^2 + y^2)|u_2|,\]
  with the boundaries $\phi_1(x,y)$ and $\phi_2(x,y)$ defined as follows:
\begin{equation*}
\phi_1(x,\pm 1)=
\begin{cases}
0  & -1 \leq x  < 0,\\
x^{1/2}   &  0 \leq x  \leq 1,
\end{cases}
\end{equation*}
 \[
 \phi_1(1,y)=1,\quad \quad \phi_1(-1,y)=0,
 \]
 and
 \begin{equation*}
\phi_2(x,\pm1)=
\begin{cases}
|x|  & -1 \leq x  < 0,\\
0   &  0 \leq x  \leq 1,
\end{cases}
\end{equation*}
\[
 \phi_2(-1,y)=1,\quad \quad \phi_2(1,y)=0.
\]
\begin{figure}[ht]
	\begin{center}
		\subfloat[$50\times 50$ disctretization]{\includegraphics[scale=.5]{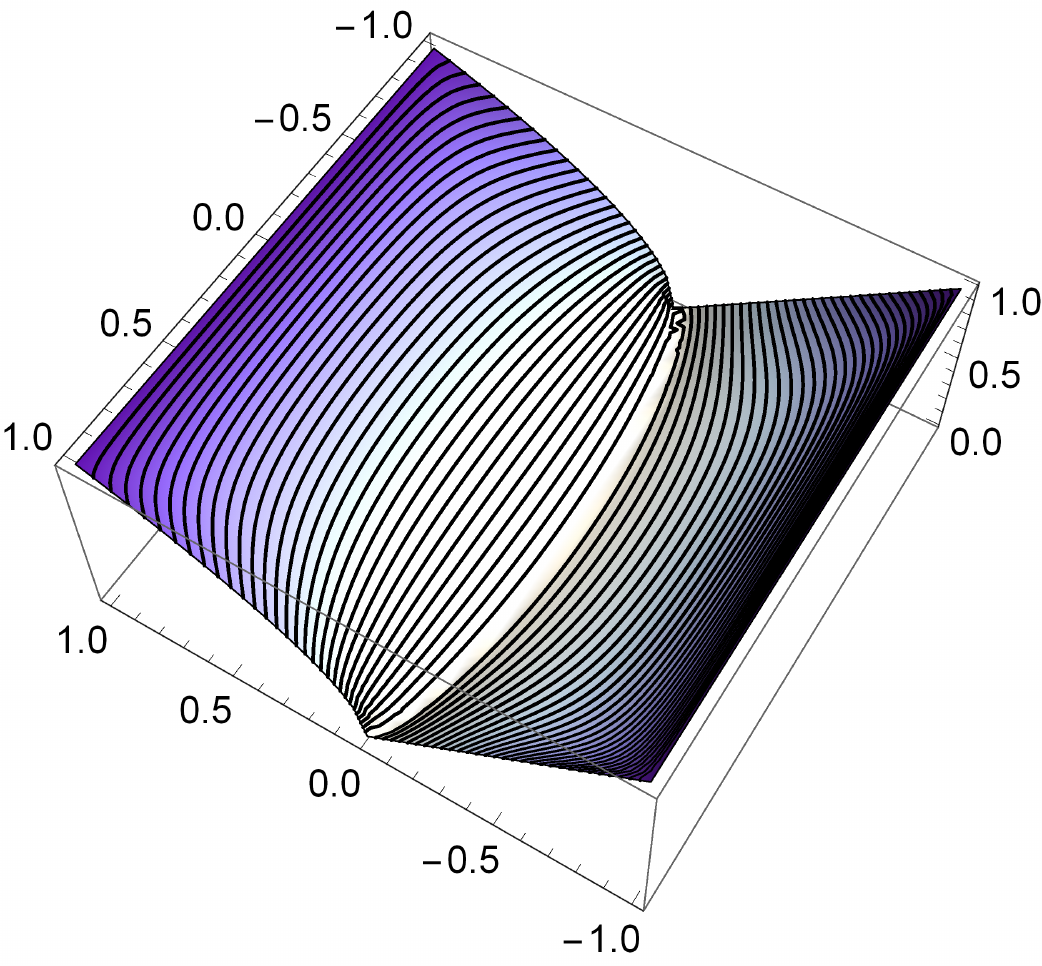}}
		\hspace{.1cm}
		\subfloat[Level sets.]{\includegraphics[scale=.5]{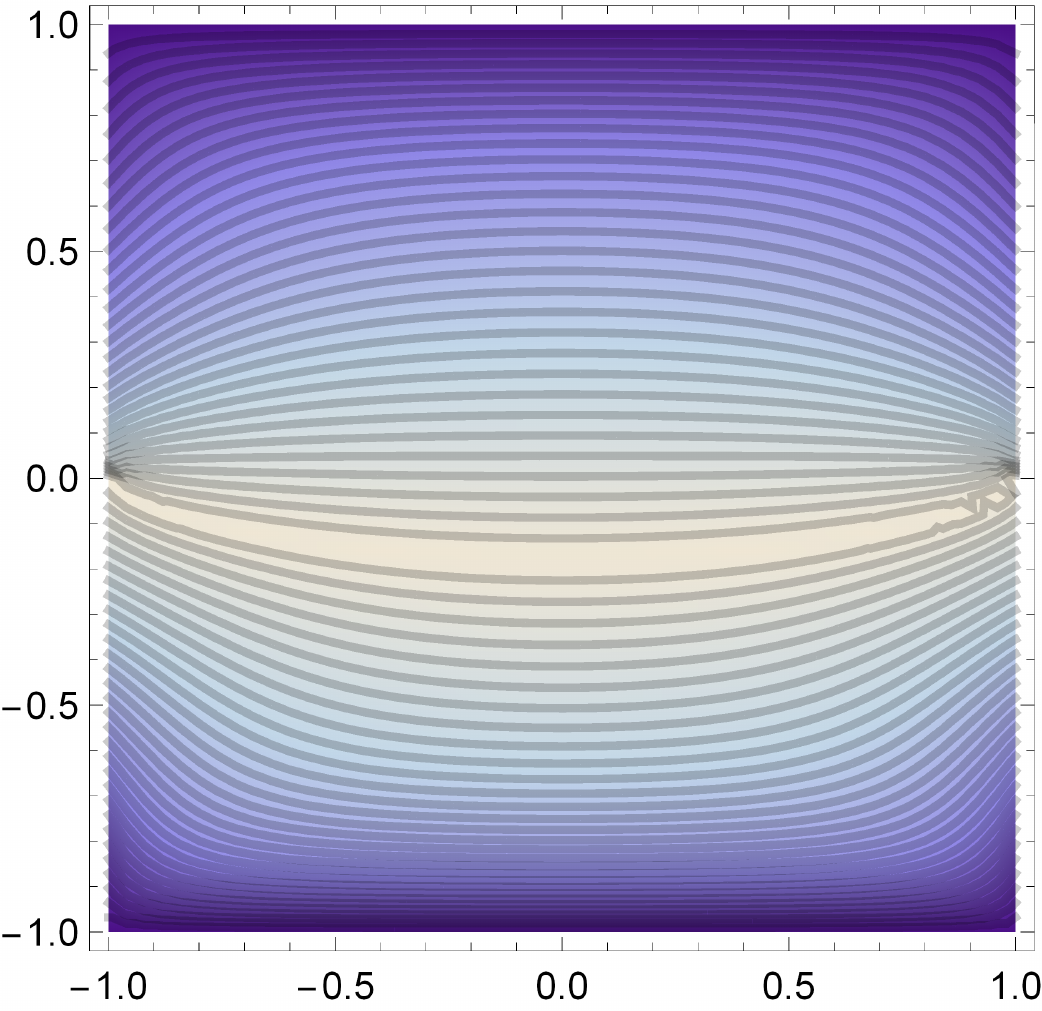}} 
	\end{center}
	\caption{Numerical solution with internal dynamics $f_1(x,y,u_1)$ and $f_2(x,y,u_2)$.}\label{twopop}
\end{figure}
In Figure \ref{twopop}, the numerical solution and its level sets are presented for $50\times 50$ discretization points and $1500$ iterations.  

Next, to show the accuracy of the iterative algorithm  in Example \ref{exampleexact}, we illustrate the difference of the exact and numerical solution.  
\end{example}

\begin{example}\label{exampleexact}
	We again take  $\Omega=[-1,1]\times[-1,1]$ and the internal dynamics $f_i(x,y,u_i)$ as follows
	  \[
	  f_1(x,y,u_1)=2,\quad f_2(x,y,u_2)=2,\quad f_3(x,y,u_3)=8,
	  \]
	  with the boundaries $\phi_1(x,y)$, $\phi_2(x,y)$ and $\phi_3(x,y)$ depicted in Figure \ref{fig1A}. Their explicit form on the boundary is the following:
\begin{equation*}
\phi_1(x,1)=
\begin{cases}
 x (x+3) & x\geq 0, \\
 0 & x < 0, \\
\end{cases}\quad\phi_1(1,x)=
\begin{cases}
 1 + 3 x & 1\geq -3x, \\
 0 & 1 <- 3x, \\
\end{cases}
\end{equation*}

\begin{equation*}
\phi_1(x,-1)=0,\quad\phi_1(-1,x)=0,
\end{equation*}

\begin{figure}[ht]
		\begin{center}
			\subfloat[$\phi_1+\phi_2+\phi_3$]{\includegraphics[scale=.5]{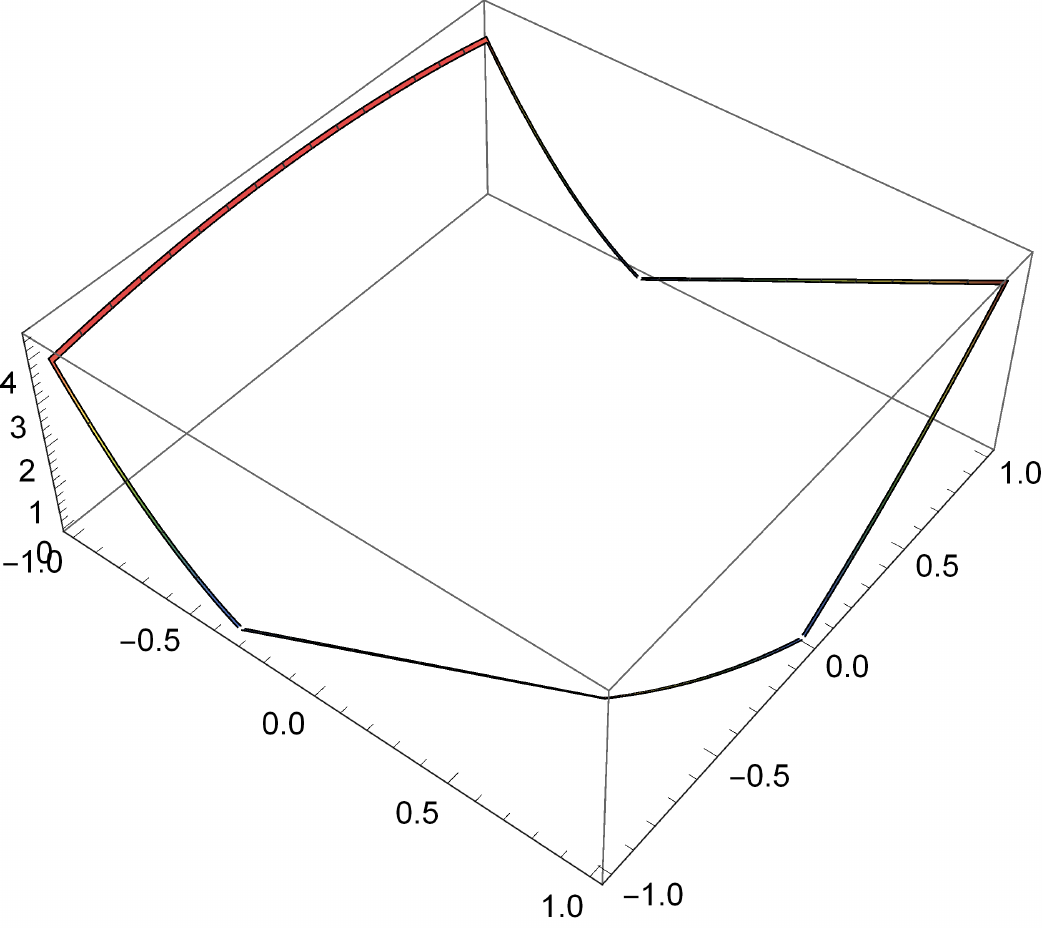}\label{fig1A}} 
			\hspace{.1cm}
			\subfloat[$u_1+u_2+u_3$]{\includegraphics[scale=.5]{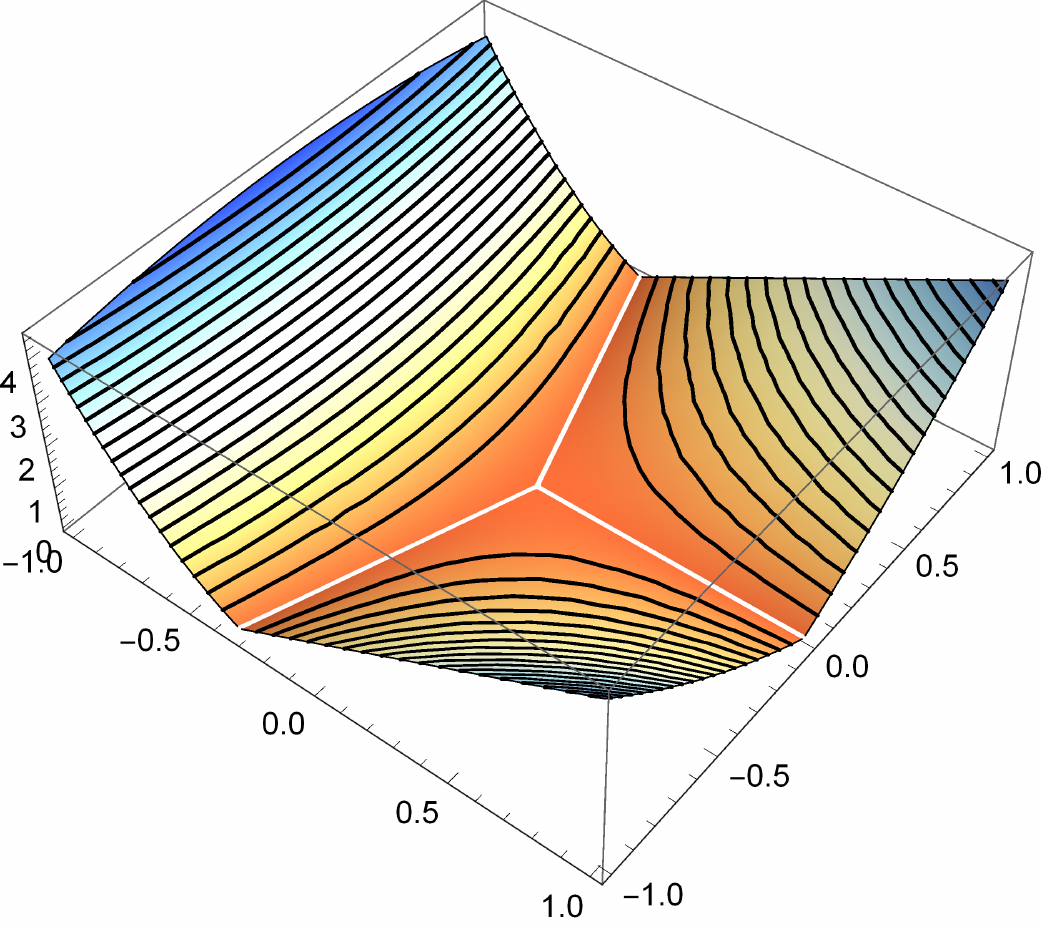}\label{fig1B}}
		\end{center}
		\caption{The sum of boundary functions $\phi_i(x,y )$ and the exact population densities $u_i(x,y)$.}\label{treepop}
\end{figure}
\begin{equation*}
\phi_2(x,1)=
\begin{cases}
(-3 + x) x & x< 0, \\
 0 & x \geq 0, \\
\end{cases}\quad\phi_2(-1,x)=
\begin{cases}
 1 + 3 x & -1\leq 3x, \\
 0 & -1> 3x, \\
\end{cases}
\end{equation*}
\begin{equation*}
\phi_2(x,-1)=0,\quad\phi_2(1,x)=0,
\end{equation*}

	\begin{equation*}
	\phi_3(x,\pm 1)=
	\begin{cases}
	\frac{1}{2} (-1 + 9 x^2) & x< -\frac{1}{3}, \\
	 0 & x\geq -\frac{1}{3}, \\
	\end{cases}
    \end{equation*}
	
	\begin{equation*}
	\phi_3(-1,x)=-\frac{1}{2} (x-3)(x+3)\quad\phi_3,\\
	(1,x)=0.
	\end{equation*}  

\begin{figure}[ht]
		\begin{center}
	       {\subfloat[$20\times 20$ disctretization ]{
			\includegraphics[scale=.5]{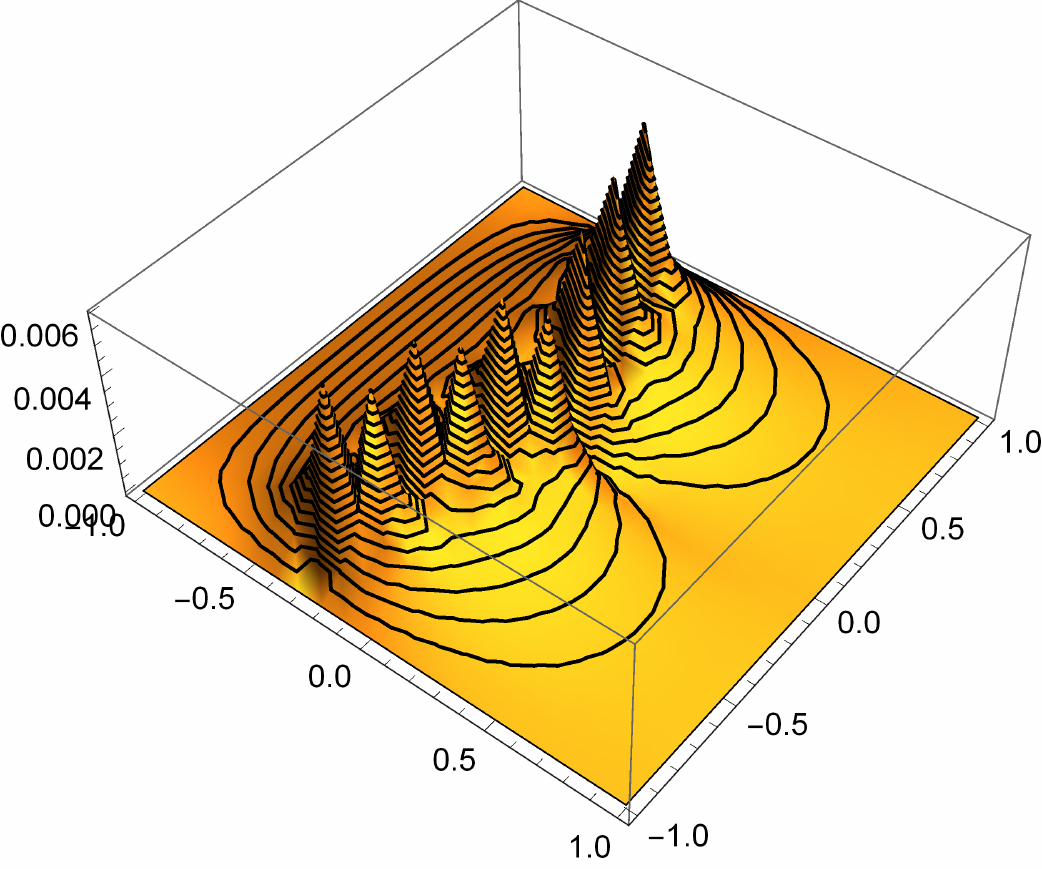}
	        \label{treepop_errorA}}
			}
	        \hspace{.1cm}
			{
	        \subfloat[$40\times 40$ disctretization]{
	        \includegraphics[scale=.5]{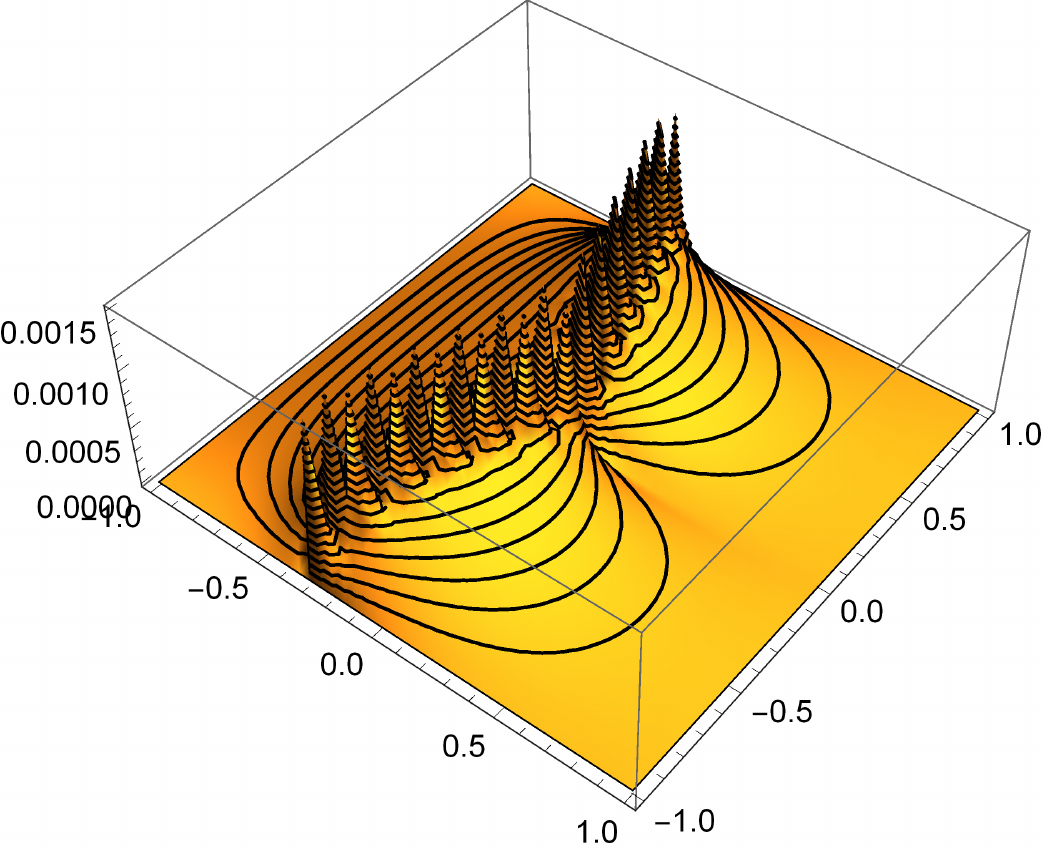}
	        \label{treepop_errorB}}
	        }
	\caption{The difference $|w-w_h|$\label{treepop_error}}
		\end{center}
\end{figure}
The exact unique minimizer $(u_1(x,y),u_2(x,y),u_3(x,y))$ to functional \eqref{last}  will be 
\begin{eqnarray*}
     u_1(x,y)=&-(3 x-y) y\quad & y< 0 \text{ and } y\leq 3 x, \\
     u_2(x,y)=&y (3 x+y)\quad & y\geq 0 \text{ and } y\geq -3 x. \\
     u_3(x,y)=&\frac{1}{2} (9 x^2 - y^2)\quad & 3 x< y \text{ and } y< -3 x,
\end{eqnarray*}
 In Figure \ref{fig1B}  it is clearly visible that the free boundaries corresponding to these populations are three straight lines emanating from the origin.

\begin{table}[ht]
\begin{center}
    \caption{Error between the exact and numerical solutions $R_{N, m\times N}$}
    \label{tb1}
    \begin{tabular}{|c|c|c|c|c|c|c|}
    \hline
    &$N=10$ & $N=20$ & $N=40$ & $N=80$ \\
    \hline
    $R_{N, 5\times N}$&$2.27 \cdot 10^{-2}$ & $6.25\cdot 10^{-3}$ & $5.06\cdot 10^{-2}$ & $3.04\cdot 10^{-1}$ \\
    \hline
    $R_{N, 10\times N}$&$2.27 \cdot 10^{-2}$ & $5.97\cdot 10^{-3}$ & $1.52\cdot 10^{-3}$ & $9.67\cdot 10^{-3}$ \\
    \hline
    $R_{N, 20\times N}$&$2.27 \cdot 10^{-2}$ & $5.97\cdot 10^{-3}$ & $1.52\cdot 10^{-3}$ & $1.64\cdot 10^{-3}$ \\
    \hline
    $R_{N, 40\times N}$&$2.27 \cdot 10^{-2}$ & $5.97\cdot 10^{-3}$ & $1.52\cdot 10^{-3}$ & $3.82\cdot 10^{-4}$ \\
    \hline
    $R_{N, 80\times N}$&$2.27 \cdot 10^{-2}$ & $5.97\cdot 10^{-3}$ & $1.52\cdot 10^{-3}$ & $3.82\cdot 10^{-4}$ \\
    \hline
    $R_{N, 160\times N}$&$2.27 \cdot 10^{-2}$ & $5.97\cdot 10^{-3}$ & $1.52\cdot 10^{-3}$ & $3.82\cdot 10^{-4}$\\
    \hline    
    \end{tabular}
\end{center}
\end{table}

We set by $w(x,y)=u_1(x,y)+u_2(x,y)+u_3(x,y),$ and by $w_h$  a numerical solution corresponding to $w$.  	In Figure \ref{treepop_error} the difference $|w-w_h|$ is presented for two various discretizations.  We consider $20\times 20$ discretization in Figure \ref{treepop_errorA} and  $40\times 40$ discretization in Figure \ref{treepop_errorB}. We clearly see that the error between  exact and numerical solutions is of the order $10^{-3},$ nevertheless its accuracy improved four times when the number of discretization points increased from $20\times 20 $ to $40\times 40$. 

Table \ref{tb1} and \ref{tb2}  show the  maximal errors between the exact and numerical solutions for  Example \ref{exampleexact}. We present the results for different discretization points  and iterations. Here  $R_{N,M}$ is the maximal error while using $N\times N$ discretization points and $M$ iterations for the algorithm \eqref{algorithm}. It is clearly visible that the error decreases with the rate around $(1/N)^2$, when the numbers $N$ and $M$ are increased (observe that when we increase the number of discretization points $2$ times, then the error decreases $4$ times).
\begin{table}[ht]
\begin{center}
    \caption{Error between the exact and numerical solutions $R_{N, m\times N}$}
    \label{tb2}
    \begin{tabular}{|c|c|c|c|c|c|c|}
    \hline
    & $N=160$& $N=320$\\
    \hline
    $R_{N, 5\times N}$& $7.79\cdot 10^{-1}$& $1.34$\\
    \hline
    $R_{N, 10\times N}$& $1.28\cdot 10^{-1}$& $4.91\cdot 10^{-1}$\\
    \hline
    $R_{N, 20\times N}$& $4.38\cdot 10^{-3}$& $8.35\cdot 10^{-2}$\\
    \hline
    $R_{N, 40\times N}$& $9.67\cdot 10^{-5}$& $3.01\cdot 10^{-3}$\\
    \hline
    $R_{N, 80\times N}$& $9.58\cdot 10^{-5}$& $2.4\cdot 10^{-5} $\\
    \hline
    $R_{N, 160\times N}$& $9.58\cdot 10^{-5}$& $2.4 \cdot 10^{-5}$\\
    \hline    
    \end{tabular}
\end{center}
\end{table}

We also see that  for both cases, the main deviation arises on two free boundary lines out of three. This is expected, since the free boundary line passing from $(0,0)$ to $(1,0)$ is falling on the net of grid points $\mathcal N$.  
\end{example}

  	\begin{figure}[ht]
  		\begin{center}
  			\subfloat[$u_1+u_2+u_3$]{\includegraphics[scale=.55]{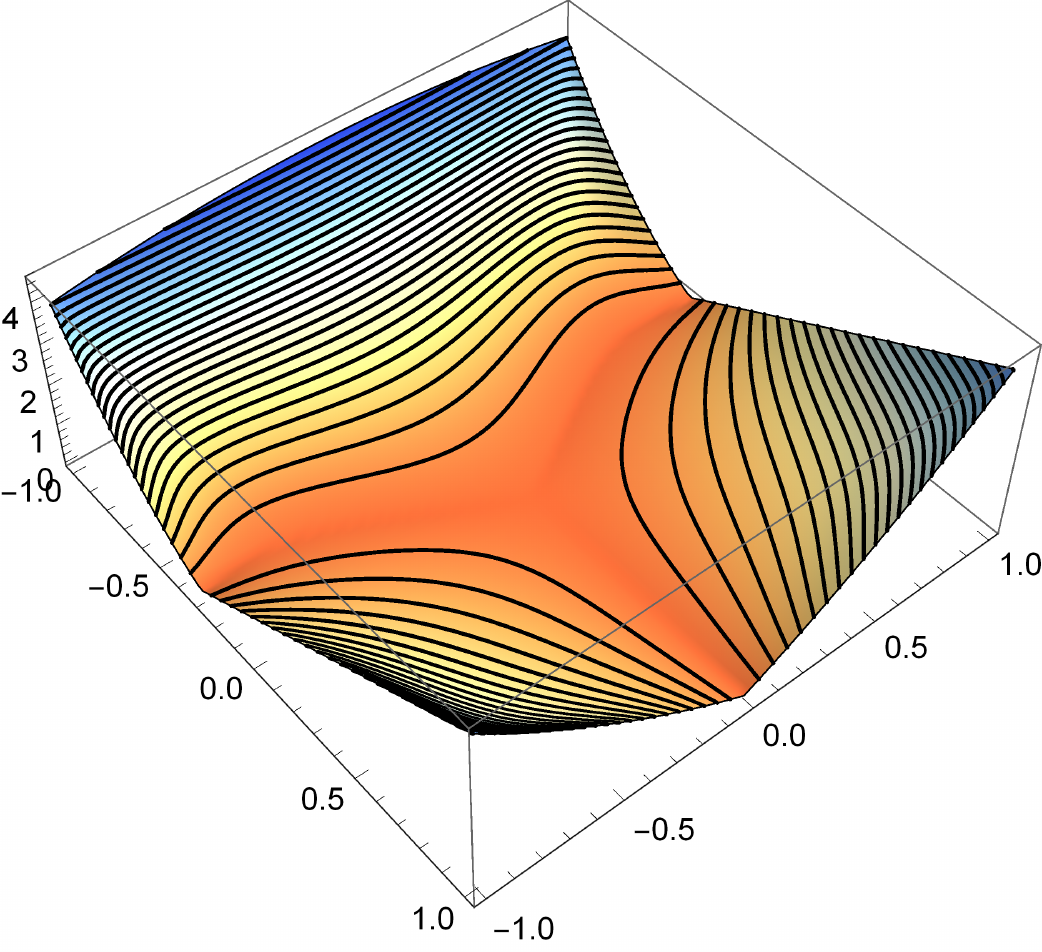}} 
  			\hspace{0.1cm}
  			\subfloat[Level sets]{\includegraphics[scale=.5]{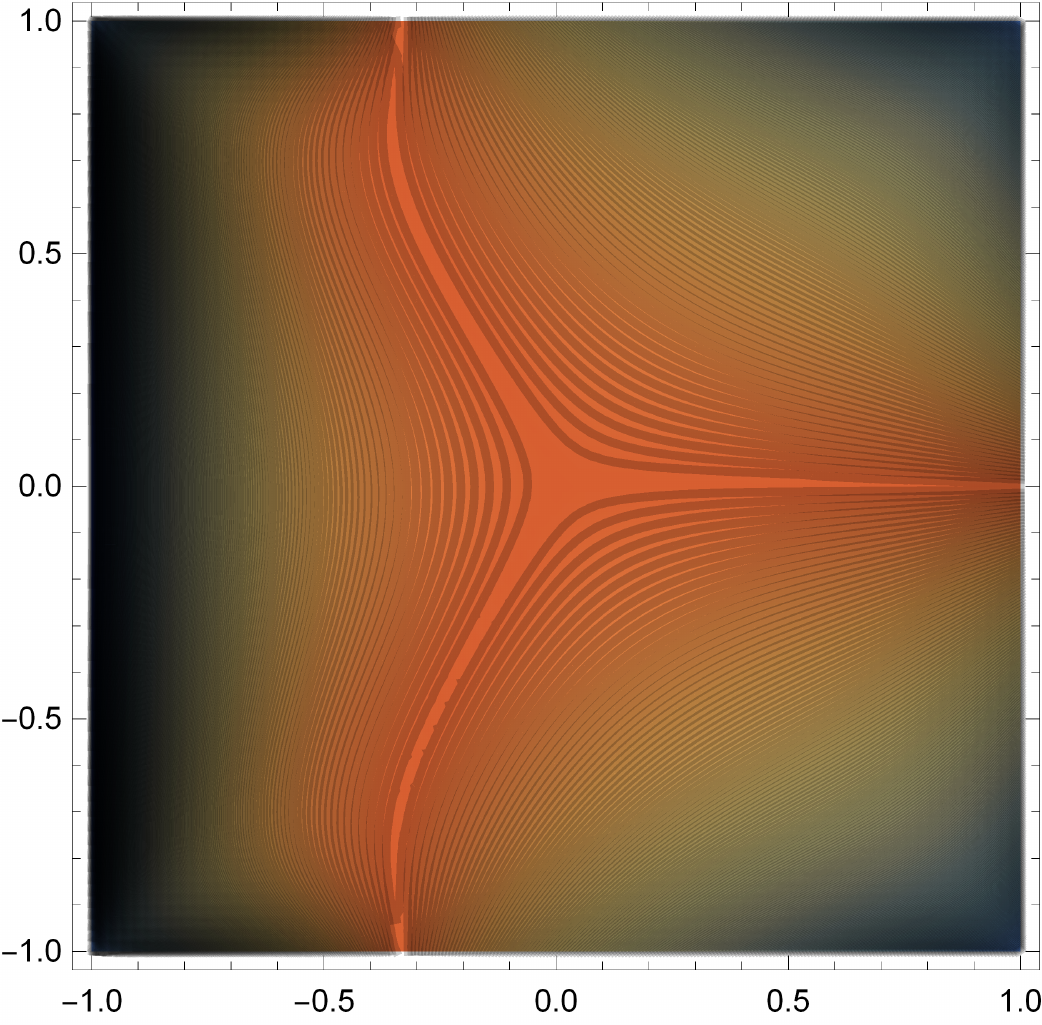}}
  		\end{center}
  		\caption{A sum of population densities $u_i(x,y)$.}\label{treepop_zero}
  	\end{figure}	
\begin{example}
	We consider three competing populations on $\Omega=[-1,1]\times[-1,1]$. The boundary values are taken as in Example \ref{exampleexact}. The difference is the internal dynamics given below
	  \[f_1(x,y,u_1)=10 (x^2 + y^2) \sqrt{u_1},\quad f_2(x,y,u_2)=10 (x^2 + y^2) \sqrt{u_2},\]
	  \[ f_3(x,y,u_3)=40 (x^2 + y^2) \sqrt{u_3}.\]

	 In this case it is hard to write down the explicit form of populations $u_i(x,y)$.  The numerical solution and its level sets are shown in Figure \ref{treepop_zero} with the use of $80\times 80$ discretization points and after $3200$ iterations. We observe that in this case, the zero set between population densities $u_i(x,y)$ is larger than in Example \ref{exampleexact}. 
\end{example}

\begin{section}{Acknowledgment}
This work was supported by State Committee of Science MES RA, in frame of the research project No.  SCS 13YR-1A0038.
The authors would like to thank L.~ Poghosyan for fruitful discussions.
\end{section}

\bibliographystyle{acm}

\bibliography{multi_obst}

\end{document}